\documentclass[11pt]{article}

\usepackage{multirow, multicol}
\usepackage{booktabs} % for much better looking tables
\usepackage{array} % for better arrays (eg matrices) in maths
\usepackage{paralist} % very flexible & customisable lists (eg. enumerate/itemize, etc.)
\usepackage{verbatim} % adds environment for commenting out blocks of text & for better verbatim
\usepackage{subfig} % make it possible to include more than one captioned figure/table in a single float
\usepackage{amsmath, amssymb}
\usepackage{pictexwd,dcpic}
\usepackage{enumerate}

\newtheorem{theorem}{Theorem}
\newtheorem{corollary}[theorem]{Corollary}
\newtheorem{remark}[theorem]{Remark}
\newtheorem{lemma}[theorem]{Lemma}
\newtheorem{definition}[theorem]{Definition}
\newtheorem{proposition}[theorem]{Proposition}
\newtheorem{example}[theorem]{Example}

\newcommand{\proof}{ {\sc Proof.\quad}}
\newcommand{\pend}{ \hfill $\square$ \\}

\numberwithin{equation}{section}  % Formeln mit fÂ¸hrender Kapitelnummer
\numberwithin{figure}{section}    % Abbildungen mit fÂ¸hrender Kapitelnummer
\numberwithin{table}{section}     % Tabellen mit fÂ¸hrender Kapitelnummer
\numberwithin{theorem}{section}

\newcommand{\of}[1]{\ensuremath{\left( #1 \right)}}

\newcommand{\cb}[1]{\ensuremath{ \left\{ #1 \right\} }}
\newcommand{\sqb}[1]{\ensuremath{ \left[ #1 \right] }}

\newcommand{\st}{\,|\;}

\newcommand{\vp}{\ensuremath{\varphi}}

\newcommand{\R}{\mathrm{I\negthinspace R}}
\newcommand{\OLR}{\overline{\mathrm{I\negthinspace R}}}
\newcommand{\N}{\mathrm{I\negthinspace N}}

\renewcommand{\P}{\ensuremath{\mathcal{P}}}

\newcommand{\G}{\ensuremath{\mathcal{G}}}

\newcommand{\C}{\ensuremath{C^-\setminus\cb{0}}}

\newcommand{\wEff}{{\rm wEff}}

\newcommand{\dom}{{\rm dom \,}}

\newcommand{\gr}{{\rm graph \,}}
\newcommand{\cl}{{\rm cl \,}}

\newcommand{\co}{{\rm co \,}}
\newcommand{\cone}{{\rm cone\,}}

\newcommand{\Int}{{\rm int\,}}

\newcommand{\isum}{{+^{\negmedspace\centerdot\,}}}

\newcommand{\idif}{{-^{\negmedspace\centerdot\,}}}

\newcommand{\lel}{\preccurlyeq}

\newcommand{\triup}{{\rm \vartriangle}}

\usepackage[%backref=page,
colorlinks=true, linkcolor=blue, citecolor=blue, urlcolor=blue, filecolor=blue
]{hyperref}

\usepackage[final]{graphicx}

\begin{document}
\title{Applying Set Optimization to Weak Efficiency}

\author{
Giovanni P. Crespi\thanks{University of Insubria, Department of Economics,  Via Monte Generoso, 71, 21100, Varese, Italy
\href{mailto:giovanni.crespi@uninsubria.it}{giovanni.crespi@uninsubria.it}}
\and
Carola Schrage\thanks{Free University of Bozen-Bolzano,
Faculty of Economics and Management, Piazza Universit\`a 1, 39100 Bozen-Bolzano Italy.
\href{mailto:carolaschrage@gmail.com}{carolaschrage@gmail.com}}%\footnote{the research was partly done while the second authors stay at the University of Valle d'Aosta}
}
\date{{\small \today}}
\maketitle
\begin{abstract}
Since the seminal papers by Giannessi, an interesting topic in vector optimization has been the characterization of (weak) efficiency through Minty and Stampacchia type variational inequalities. Several results have been proved to extend those known for the scalar case. However, in order to introduce a proper definition of variational inequality, some assumptions are usually made that may eventually be questioned.\\
We find two major approaches  in the papers we considered, that arise when defining generalized derivatives for vector--valued functions. First, some authors introduce set--valued derivatives for single--valued problems, thus completely changing the setting of the problem.

Second, when dealing with Dini--type derivatives, infinite elements may occur. The approach to handle this problem is not uniquely defined in the literature, therefore, when considered, the definition proposed may seem arbitrary.\\
Indeed these problems are strictly related to the lack of a complete order in the image space of a vector--valued function. We use an alternative approach to study vector optimization, by considering a set--valued counterpart defined with values in a conlinear space. The structure of this space allows to overcome the previous difficulties and to obtain variational inequality characterization of weak efficiency as a straightforward application of scalar arguments.
\end{abstract}

%%%%%%%%%%%%%%%%%%%%%%%%%%%%%%%%%%%%%%%%%%%%%

\section{Introduction}\label{sec:Intro}

%Vector optimization has been extensively studied in the literature.
	Since the seminal papers by Giannessi \cite{Giannessi80, Giannessi97} one of the issues in (convex) vector optimization has been the use of \emph{differentiable variational inequalities} to characterize weak efficient solutions of a \emph{primitive} optimization problem, see e.g. \cite{CrespiGinRoc2005, Ginchev2007}. Given a differentiable, vector valued, objective function $\psi:S\subseteq X\to Z$, where $X$ and $Z$ are vector spaces, and the partial order induced by a closed, convex, pointed cone with nonempty interior $C\subset Z$, the vector optimization problem is 
		\begin{align}\label{VOP}\tag{VOP}
			\min\psi(x),\; x\in S
		\end{align}
	A weak efficient solution \eqref{VOP} is $x_0\in S$ such that $\psi(x_0)$ is a weakly efficient element in the image set $\psi\sqb{X}=\cb{\psi(x)\st x\in S}$, i.e.
		\begin{align*}
			\forall z\in \psi\sqb{X}:\quad \psi(x_0)\notin \cb{z}+\Int C.
		\end{align*}
	Problem \eqref{VOP} is referred to as \emph{primitive} when compared to the Minty variational inequality problem of finding $x_0\in S$ s.t.
		\begin{align}\label{MVIP}\tag{MVIP}
			\langle\psi^\prime(x),x-x_0\rangle\in-\Int C\;\forall x\in S
		\end{align}
	or the Stampacchia variational inequality problem of finding $x_0\in S$ s.t.
		\begin{align}\label{SVIP}\tag{SVIP}
			\langle\psi^\prime(x_0),x_0-x\rangle\in-\Int C\;\forall x\in S
		\end{align}
	Since the variational inequalities \eqref{MVIP} and \eqref{SVIP} are defined through the derivative $\psi^\prime$, they are usually referred to as \emph{differentiable}. Indeed, it is rather obvious that both variational inequalities define directional derivatives of the differentiable objective function (hence \emph{primitive} to the inequality). Therefore, the problem has been soon after extended to the nondifferentiable case by using generalized directional derivatives to replace the inner products in the previous formulation. Relations between the set of weak efficient solutions of \eqref{VOP} and those of the associated (generalized) variational inequalities have been proved in various papers compare e.g. \cite{AnsariLee2010, cgr-erice, YYT}.
	The goal to be achieved in vector optimization is to extend known results for scalar optimization, where, starting from the convex, differentiable case, more general results have been achieved. In \cite{cgrJOTA}, for instance, under mild assumptions, solutions to a Minty-type variational inequality are proven to be global minimizers of the underlaying non-convex optimization problem. However, when vector optimization is involved, several instances has to be considered. Among others, the directional derivatives that has been proved effective in the scalar case to tackle the non differentiable case, involve upper or lower limits of difference quotients, that are not straightforward for vector-valued functions, mainly due to the lack of completeness in the order induced by $C$. To overcome these difficulties, we have found two major approaches in the literature. One involves the definition of set--valued derivative for single--valued problems (see e.g. \cite{Ginchev2007}), thus completely changing the image space of \eqref{VOP}. The other, to copy with the possibility of non finite generalized directional derivatives, requires to introduce arbitrary notions of infinite elements for vector spaces (e.g. \cite{Ginchev2007}), or to simply avoided the possibility, imposing finiteness of the limits (e.g. \cite{AnsariLee2010}).
	
	In this paper we restrict ourselves to the convex case, in order to better exploit the problem of differentiability. It is left as an open question and future line of research the non convex case and the possible generalizations of results in \cite{cgrJOTA}. Along the lines of \cite{HeydeLoehne11}, where a 'fresh look' to vector optimization has been proposed by means of set-optimization, we propose a suitable set-optimization problem to study weak efficiency in \eqref{VOP}. This approach allows to overcome the ambiguity of infinite elements, dealing with a fully set--valued problem, and gaining a deeper insight on the original vector--valued problem. 
	To gain such a result, we side $\psi$ with its set--valued extension, $\psi^C$, mapping $S$ into an order complete space $\G^\triup$, as in \cite{galatos2007residuated, GetanMaLeSi, HamelHabil05, HamelSchrage12, Loehne11Book, MartinezLegazSinger95}. Moreover, we define directional derivatives, mapping into $\G^\triup$, by means of upper or lower limits of difference quotients in the image space, that can be applied to $\psi^C$ in order to prove necessary and sufficient conditions in terms of Stampacchia and Minty variational inequalities, to characterize the set of weak solutions to \eqref{VOP}, under convexity assumptions.
	Eventually, results as proved in \cite{CGR, Giannessi97} follow, as a special cases, overcoming the necessity to introduce infinite elements in $Z$ and to study the topology of the extended vector space $\tilde{Z}$, that appears in the cited papers.
	Corresponding results for stronger variational inequalities and minimizer rather then weak minimizers have been studied in \cite{CrespiHamelSchrage2015} and \cite{CrespiSchrage13a}.

The remainder of the paper is organized as follows.
In Section  \ref{sec:Setting}, we introduce the general setting and the basic notation. Some results for the set $\G^\triup$ are proven for subsequent reference. Section \ref{sec:DirDer} is devoted to the concept of upper and lower Dini directional derivatives for functions mapping onto $\G^\triup$. We show that these concepts generalize the original definition for proper scalar functions, compare e.g. \cite{giorgi1992dini}.
The final Section \ref{sec:Main_Results} collects our main results, applying the general scheme to \eqref{VOP}. In this final section, we restrict ourselves to the case of convex functions in order to achieve a greater simplicity of the arguments, rather then greatest possible generality, leaving the more general case for further research.

%%%%%%%%%%%%%%%%%%%%%%%%%%%%%%%%%%%%%%%%%%%%%%%%%%%%

\section{Setting}\label{sec:Setting}

%%%%%%%%%%%%%%%%%%%%%%%%%%%
In this paper we consider locally convex Hausdorff spaces $X$ and $Z$,  with topological dual $Z^*$, and $\P(Z)$ the power set of $Z$, including $\emptyset$ and $Z$ as elements.
Throughout the paper we denote by $\mathcal U_X$ and $\mathcal U_Z$ the set of all closed, convex and balanced $0$ neighbourhoods in $X$ and $Z$ respectively and by $\cl A$, $\co A$ and $\Int A$, the closed hull, the convex hull and the topological interior of a set $A$, respectively. The conical hull of a set $A$ is $\cone A=\cb{ta\st a\in A,\, 0<t}$. To define a solution concept to (\ref{VOP}) we introduce a preorder on $Z$ by a closed convex cone $C\neq Z$ with nonempty topological interior, $\Int C\neq\emptyset$. As usual, by $z_1\leq z_2$ we mean $z_2\in \cb{z_1}+C$. The (negative) dual cone of $C$ is the set $C^-=\cb{z^*\in Z^*\st \forall z\in C:\, z^*(z)\leq 0}$. Since $\Int C\neq \emptyset$ is assumed, there exists a weak$^*$ compact base $W^*$ of $C^-$, i.e. a convex subset with $\C=\cone W^*$ with $z^*, tz^*\in W^*$ implying $t=1$ and any net in $W^*$ has a weak$^*$ convergent subnet, compare \cite[Theorem 1.5.1]{Aubin71}. Also, for every $z\in Z$ it holds $\inf\cb{w^*(z)\st w^*\in W^*}>-\infty$ and for any $U\in \mathcal U_Z$ it holds $\sup\cb{\inf\cb{w^*(u)\st u\in U}\st w^*\in W^*}<0$, compare \cite[Remark 3.32]{HeydeSchrage11R}.\\
Recall that an ordering cone is Daniell, or has the Daniell property, if and only if every decreasing net which is bounded from below converges to its infimum.
A convex polyhedron is the intersection of finitely many closed halfspaces. Especially, a closed convex cone $C$ is polyhedral, if and only if there is a finite set $M^*\subseteq W^*$ such that $C=\bigcap_{m^*\in M^*}\cb{z\in Z\st m^*(z)\leq 0}$.%, in which case $W^*$ is a polytope.
In the sequel, given any vector--valued function $\psi:S\subseteq X\to Z$, we define its set--valued extension $\psi^C:X\to\P(Z)$ as the function mapping $x$ to the upper Dedekind cut of $\psi(x)$ with respect to $C$, namely
$$
\psi^C(x)=\left\{
\begin{array}{ll}
\cb{\psi(x)}+C & \mbox{\rm if\;} x\in S\\
\emptyset & \mbox{\rm elsewhere.}
\end{array}\right.
$$

Images of $\psi^C:X\to \P(Z)$ are closed convex sets, closed under the addition with the ordering cone $C$. %, that is $\psi^C(x)=\cl\co(\psi^C(x)+C)$.
Therefore we restrict our focus to the set
\[
\G(Z,C)=\cb{A\in\P(Z)| A=\cl\co(A+C)}
\]
as a natural image space for the set-valued functions throughout this paper. Properties of $\G(Z,C)$ have been extensively studied in recent years,
compare \cite{ HamelHabil05, FivePersonSoMF, Loehne11Book, Diss}.
First we recall that the ordering in $Z$ can be extended to the power set of $Z$ (compare \cite{Hamel09, Kuroiwa98-1}
and the references therein) by setting
\[
A_1\lel A_2\quad\Leftrightarrow\quad A_2+C\subseteq A_1+C
\]
for all $A_1, A_2\subseteq Z$.  By the order relation $\lel$ defined through $\supseteq$ on the subset $\G(Z,C)$ the complete lattice $(\G(Z,C),\supseteq)$ has been introduced, compare e.g. \cite{HamelSchrage12}.
	%Especially, the relation $\lel$ coincides with $\supseteq$ on the subset $\G(Z,C)$ and $(\G(Z,C),\supseteq)$ is a complete lattice, see e.g. \cite{HamelSchrage12}.

In this framework, it is an easy exercise to apply basic set theory to the given setting to prove that for any subset $\mathcal A\subseteq \G(Z,C)$, supremum and infimum of $\mathcal A$ in $\G(Z,C)$ are given by
\begin{align*}
\inf\mathcal A=\cl\co \bigcup\limits_{A\in \mathcal A}A;\quad \sup\mathcal A=\bigcap\limits_{A\in \mathcal A}A,
\end{align*}
compare \cite[Proposition 5.18]{HamelSchrage12}.

When $\mathcal A=\emptyset$, we agree that $\inf\mathcal A=\emptyset$ and $\sup\mathcal A =Z$. Hence $\G(Z,C)$ possesses a greatest and smallest element $\inf\G(Z,C)=Z$ and $\sup\G(Z,C)=\emptyset$.
The Minkowsky addition and multiplication with non--negative reals need to be slightly adjusted to provide operations on $\G(Z,C)$. We define
\begin{align*}
\forall A, B\in \G(Z,C):\quad 		&A\oplus B=\cl\cb{a+b\in Z\st a \in A,\, b\in B};\\
\forall A\in\G(Z,C),\, \forall 0<t:\quad &t\cdot A=\cb{ta \in Z\st a\in A };\quad 0\cdot A=C.
\end{align*}
Note that $0\cdot \emptyset=0\cdot Z=C$ and $\emptyset$ dominates the addition in the sense that $A\oplus \emptyset=\emptyset$ is true for all $A\in\G(Z,C)$.
Moreover, $A\oplus C=A$ is satisfied for all $A\in \G(Z,C)$, thus $C$ is the neutral element with respect to addition.

As a consequence,
\begin{align*}
\forall \mathcal A\subseteq \G(Z,C),\, \forall B\in\G(Z,C) :\quad B\oplus \inf\mathcal A=\inf\cb{B\oplus A\st A\in \mathcal A},
\end{align*}
or, equivalently,  the $\inf$--residual
\[
A\idif B=\inf\cb{M\in\G(Z,C)\st B\oplus M\subseteq A}
\]
exists for all $A, B\in \G(Z,C)$. It holds (compare \cite[Theorem 2.1]{HamelSchrage12})
\begin{align*}
A\idif B    &=\cb{z\in Z\st B+\cb{z}\subseteq A};\\
A             &\supseteq B\oplus (A\idif B).
\end{align*}

Overall, the structure of $\G^\triup=\of{\G(Z,C),\oplus,\cdot,C,\supseteq}$ is that of an order complete $\inf$--residuated conlinear space.
Since the seminal paper \cite{HamelHabil05}, conlinear spaces have been throughly studied. Residuation is well known in order theory, compare \cite{Fuchs66, galatos2007residuated} and has been applied to convex analysis by Mart\'inez-Legaz, Singer and Getan, compare \cite{GetanMaLeSi, MartinezLegazSinger95}. We further remark that residuation provides a substitute for the difference operation and becomes a powerful tool to extend calculus to set-valued functions.
For the reader convenience, we briefly recall the definition.
%{\color{red}
\begin{definition}
\label{DefConlinearSpace} A nonempty set $Y$ together with two algebraic operations $+
: Y \times Y \to Y$ and $\cdot : \R_+ \times Y \to Y$  is called a conlinear
space with neutral element $\theta$ provided that
\\
(C1) $\of{Y, +,\theta}$ is a commutative monoid with neutral element $\theta$: For all $ w_1, w_2, w_3 \in Y$ it holds (i) $w_1+ w_2=w_2+ w_1\in Y$, (ii) $w_1+(w_2+w_3)=(w_1+w_2)+w_3$, (iii)$w_1+\theta=\theta+w_1=w_1$;
\\
(C2) The operations are compatible: (i) $\forall w_1, w_2 \in Y$, $\forall r \in \R_+$: $r \cdot \of{w_1 + w_2} = r
\cdot w_1 + r \cdot w_2$, (ii) $\forall w \in Y$, $\forall r, s \in \R_+$: $s \cdot \of{r
\cdot w} = \of{rs} \cdot w$, (iii) $\forall w \in Y$: $1 \cdot w = w$, (iv)  $\forall w \in Y$: $0 \cdot
w = \theta$.
\\
Subsequently, these operations are referred to as addition and multiplication, respectively.
\\
A conlinear space $\of{Y, +, \cdot,\theta}$ together with an order relation $\lel$ on $Y$ is called partially ordered, lattice ordered or order complete  conlinear space provided
that $(Y,\lel)$ has the respective structure and  the order is compatible with addition and multiplication, that is\\
(C3) (i) $\forall w, w_1, w_2 \in Y$, $w_1 \lel w_2$ imply $w_1 + w \lel w_2 + w$,
 and (ii) $\forall w_1, w_2 \in Y$, $w_1 \lel w_2$, $r \in \R_+$ imply $r \cdot w_1 \lel r\cdot w_2$.

A partially ordered conlinear space  $\of{Y, +, \cdot,\theta,\lel }$ is called $\inf$--residuated, when for all $w_1,w_2\in Y$ the element
$w_2\idif w_1=\inf\cb{u\in Y\st w_2\lel w_1+ u}$
exists. In this case, $w_2\idif w_1$ is called the $\inf$--residual of $w_2$ and $w_1$.
\end{definition}

A partially ordered conlinear space $Y$ is $\inf$--residuated, if and only if for all $w\in Y$ and all $A\subseteq Y$ such that $\inf  A$ exists, it holds $\of{w+ \inf A}=\inf\cb{w+ a\st a\in A}$ (compare \cite[Theorem 2.1]{HamelSchrage12}).
%}

References and details on  structural properties of conlinear spaces and $\inf$--residuation can be found
in \cite{galatos2007residuated, HamelHabil05, FivePersonSoMF, HamelSchrage12}.

%{\color{red}
Notably,
$(s+t)x=sx+tx$ for $s,t\in\R_+$ and $x\in Y$ is not assumed on a conlinear space and is not satisfied for $\G(Z,C)$. Thus, the power set of a conlinear space is again a conlinear space.  The lack of the second associativity law is what sets conlinear spaces apart from other concepts, such as semilinear spaces \cite[p.145]{kutateladze1972minkowski} or abstract convex cones \cite{FuchssteinerLusky1981}.
%}

The following result can easily be extended to general $\inf$--residuated conlinear spaces.

\begin{lemma}\label{lem:calc_conlin}
Let $A,B,D,E\in \G^\triup$ and $0<t,s\in\R$. Then
\[
(tA\oplus sB)\idif (tD\oplus sE)\supseteq t(A\idif D)\oplus s(B\idif E).
\]
\end{lemma}
\proof
Since the ordering in $\G^\triup$ is compatible with the algebraic operations and $t(A\idif B)=tA\idif tB$ is true for all $0<t$, without loss of generality we can assume $t=s=1$.
As $\G^\triup$ is $\inf$--residuated,
\begin{align*}
(A\idif D)\oplus (B\idif E)&=\inf\cb{T\in \G^\triup\st A\supseteq D\oplus T}\oplus\inf\cb{S\in \G^\triup\st B\supseteq E\oplus S}\\
&=\inf\cb{T\oplus\inf\cb{S\in \G^\triup\st B\supseteq E\oplus S}\in \G^\triup\st A\supseteq D\oplus T}\\
&=\inf\cb{T\oplus S \in \G^\triup\st A\supseteq D\oplus T,\, B\supseteq E\oplus S}
\end{align*}
but $A\supseteq D\oplus T$ and $B\supseteq E\oplus S$ together imply
\[
A\oplus B\supseteq (D\oplus T)\oplus(E\oplus S)=(D\oplus E)\oplus(T\oplus S),
\]
hence
\[
A\oplus B\supseteq (D\oplus E)\oplus((A\idif D)\oplus (B\idif D)),
\]
and equivalently
\[
(A\oplus B)\idif (D\oplus E)\supseteq (A\idif D)\oplus (B\idif E).
\]
\pend

In the sequel, we will make use of the fact that $\OLR$ equipped with adequate ordering and addition can be identified with the space $\G(\R,\R_+)$.

\begin{example}
\label{ExExtReals}
Let us consider $Z = \R$, $C = \R_+$. Then $\G\of{Z, C} = \cb{[r, +\infty) \mid r \in \R}\cup\cb{\R}\cup\cb{\emptyset}$, and $\G^\triup$ can be identified (with respect to the algebraic and order structures which turn $\G\of{\R, \R_+}$ into an ordered conlinear space and a  complete lattice admitting an inf-residuation) with $\OLR = \R\cup\cb{\pm\infty}$ using the 'inf-addition' $\isum$ (see \cite{HamelSchrage12})%RockafellarWets98
\[
r \isum s  = \inf\cb{a+b\in\R \mid r \leq a,\, s\leq b}.
\]
The inf-residuation on $\OLR$ is given by
\[
r \idif s  = \inf\cb{t\in\R \mid r \leq s \isum t}
\]
for all $r,s\in\OLR$, compare \cite{HamelSchrage12} for further details.
\end{example}

%{\color{red}
 Since each element of $\G^\triup$ is closed and convex and $A=A+C$, by a separation argument $A$ is equal to the closed halfspaces containing it, hence
% }
\begin{align}\label{eq:scal_representation_Set}
\forall A\in \G^\triup:\quad A=\bigcap\limits_{z^*\in W^*}\cb{z\in Z\st -\sigma(z^*| A)\leq -z^*(z)},
\end{align}
where $\sigma(z^*| A)=\sup\cb{z^*(z)\st z\in A}$ is the support function of $A$ at $z^*$.
%\begin{remark}\rm
%$A=\emptyset$ if and only if there exists $z^*\in W^*$ such that $-\sigma(z^*| A)=+\infty$, or equivalently if the same holds true for all $z^*\in W^*$.
%\end{remark}

%The right hand side in \eqref{eq:scal_representation_Set} provides a scalarization of elements $A\in\G^\triup$.
The following equivalent formulation holds as well
\begin{align}\label{eq:RecCone_finiteScal}
\forall A\in \G^\triup\setminus\cb{\emptyset}:\quad A=\bigcap\limits_{\substack{z^*\in W^*,\\ -\sigma(z^*|A)\in\R}}\cb{z\in Z\st -\sigma(z^*| A)\leq -z^*(z)},
\end{align}

Applying these characterizations, scalarized counterparts of infimum and supremum of a subset of elements in $\G^\triup$ are provided.

\begin{lemma}\label{prop:scal_of_infimum}\cite[Proposition 3.5(c)]{Schrage10Opt}
Let $\mathcal A\subseteq\G^\triup$ be a set, then
\begin{align*}
&\inf\mathcal A=\bigcap\limits_{z^*\in W^*}\cb{z\in Z\st \inf\cb{-\sigma(z^*| A)\st A\in\mathcal A}\leq -z^*(z)}\\
&\forall z^*\in W^*:\quad -\sigma(z^*| \inf\mathcal A)=\inf\cb{-\sigma(z^*| A)\st A\in\mathcal A}.
 \end{align*}
\end{lemma}

\begin{lemma}\label{prop:scal_of_supremum}\cite[Lemma 4.14]{FivePersonSoMF}
Let $\mathcal A\subseteq\G^\triup$ be a set, then
\begin{align*}
&\sup\mathcal A=\bigcap\limits_{z^*\in W^*}\cb{z\in Z\st \sup\cb{-\sigma(z^*| A)\st A\in\mathcal A}\leq -z^*(z)}\\
&\forall z^*\in W^*:\quad -\sigma(z^*| \sup\mathcal A)\geq\sup\cb{-\sigma(z^*| A)\st A\in\mathcal A}.
 \end{align*}
\end{lemma}
%\proof
%For any set $\mathcal A\subseteq\G^\triup$, the supremum of $\mathcal A$ is given by $\sup\mathcal A=\bigcap\limits_{A\in \mathcal A}A$, and by the scalarization formula \eqref{eq:scal_representation_Set}, it holds
%\begin{align*}
%\sup\mathcal A
%&= \bigcap\limits_{A\in \mathcal A}\bigcap\limits_{z^*\in W^*}\cb{z\in Z\st -\sigma(z^*| A)\leq -z^*(z)}\\
%&= \bigcap\limits_{z^*\in W^*}\bigcap\limits_{A\in \mathcal A}\cb{z\in Z\st -\sigma(z^*| A)\leq -z^*(z)}\\
%&= \bigcap\limits_{z^*\in W^*}\cb{z\in Z\st \sup\cb{-\sigma(z^*| A)\st A\in \mathcal A}\leq -z^*(z)}.
%\end{align*}
%This implies
%\[
%\cb{z\in Z\st \sup\cb{-\sigma(z^*| A)\st A\in \mathcal A}\leq -z^*(z)}\lel \sup\mathcal A
%\]
%for all $z^*\in W^*$, hence $ -\sigma(z^*| \sup\mathcal A)\geq\sup\cb{-\sigma(z^*| A)\st A\in\mathcal A}$ is proven.
%\pend

Setting $A_i=\cb{(x_1,x_2)\in\R^2\st x_1 >  0,\, x_2\geq(i+\frac{1}{x_1})}\in\G(\R^2,\R^2_+)$ and $\mathcal A=\cb{A_i\st i\in\N}$, then $\sup\mathcal A=\emptyset$ and it is easy to see that typically the inequality in Lemma \ref{prop:scal_of_supremum} is strict, as $z^*=(-1,0)^T\in\R^2_+$ and
\[
\sup\limits_{i\in\N}-\sigma(z^*|A_i)=0 < +\infty=-\sigma(z^*|\sup\mathcal A).
\]

\begin{lemma}\label{lem:scal_of_difference}\cite[Proposition 5.20]{HamelSchrage12}
Let $A, B\in\G^\triup$, then
\begin{align*}
&A\idif B=\bigcap\limits_{z^*\in W^*}\cb{z\in Z\st (-\sigma(z^*| A))\idif (-\sigma(z^*| B))\leq -z^*(z)};\\
&\forall z^*\in W^*:\quad -\sigma(z^*| A\idif B)\geq\of{-\sigma(z^*| A)}\idif\of{-\sigma{z^*|B}}.
\end{align*}
\end{lemma}

%The recession cone of a nonempty closed convex set $A\subseteq Z$ is the closed convex cone
%$0^+A=\cb{z\in Z\st A+\cb{z}\subseteq A}$,
%compare \cite[p.6]{Zalinescu02}. By definition, $0^+\emptyset=\emptyset$ is assumed. If $A\in\G^\triup\setminus\cb{\emptyset}$, then $0^+A=A\idif A$ and $C\subseteq 0^+A$ are satisfied. Especially, $\Int\of{0^+A}\neq\emptyset$ and $(0^+A)^-\subseteq C^-$, hence $W^*\cap (0^+A)^-$ is a weak$^*$ compact base of $(0^+A)^-$.
%
%
%\begin{remark}\label{rem:Diff_of_epigraphical_ext}\rm
%Let $A, B\in\G^\triup$ be given with $A=\cb{a}+C$, $a\in Z$. Then  $0^+A=C$  and $-\sigma(z^*| A)=-z^*(a)$ is satisfied for all $z^*\in W^*$. Moreover, $B\idif A=B+\cb{-a}$ is true, hence
%\[
%\forall z^*\in W^*:\quad -\sigma(z^*| B\idif A)=\of{-\sigma(z^*| B)}\idif\of{-\sigma\of{z^*|A}}.
%\]
%\end{remark}
%{\color{red}
The recession cone of a closed convex set $A\subseteq Z$ is the set of all directions of recession of $A$,
\[
0^+A=\cb{z\in Z\st A+\cone\cb{z}\subseteq A}.
\]
The recession cone $0^+A$ is a closed convex cone, compare \cite[p.6]{Zalinescu02}.
%}
By definition we set $0^+\emptyset=\emptyset$. If $A\in\G^\triup\setminus\cb{\emptyset}$, then $0^+A=A\idif A$ and $C\subseteq 0^+A$ are satisfied. Especially, $\Int\of{0^+A}\neq\emptyset$ and $(0^+A)^-\subseteq C^-$, hence $W^*\cap (0^+A)^-$ is a weak$^*$ compact base of $(0^+A)^-$.
The recession cone of a nonempty set $A\in\G^\triup$ is what is referred to as a a 'generalized zero' in \cite{HamelSchrage13PJO}; it is the neutral element in $\G(Z,0^+A)$ and in the formulations of the variational inequalities such generalized zeros, the recession cone of images of the primitive function, will serve to replace the zero in the scalar formulation of the corresponding inequality, compare \eqref{star} and \eqref{eq:WMVI} below.

\begin{remark}\label{rem:Diff_of_epigraphical_ext}\rm
Let $A, B\in\G^\triup$ be given with $A=\cb{a}+C$, $a\in Z$. Then  $0^+A=C$  and $-\sigma(z^*| A)=-z^*(a)$ is satisfied for all $z^*\in W^*$. Moreover, $B\idif A=B+\cb{-a}$ is true, hence
\[
\forall z^*\in W^*:\quad -\sigma(z^*| B\idif A)=\of{-\sigma(z^*| B)}\idif\of{-\sigma\of{z^*|A}}.
\]
\end{remark}

The recession cone $0^+A$ of any element $A\in \G^\triup$ is related to the values of the support function of $A$ as the following two lemmas show.
\begin{lemma}\label{lem:Rec_A}
Let $A\in\G^\triup$ be a nonempty set, then
%\begin{equation}\label{eq:RecCone_finiteScal_1}
%0^+A=\cb{z\in Z\st \forall z^*\in W^*:\; -\sigma(z^*|A)=-\infty\,\vee\, 0\leq -z^*(z)}.
%\end{equation}
%Especially, for all $A\in\G^\triup$,
%either $A=\emptyset$, or
\begin{equation}\label{eq:RecCone_finiteScal_2}
0^+ A=\bigcap\limits_{\substack{z^*\in W^*\\ -\sigma(z^*| A)\in\R}}\cb{z\in Z\st 0\leq -z^*(z)}.
\end{equation}
\end{lemma}
\proof
Assume $z\notin 0^+A$, then either $A=\emptyset$ or there exists a $z^*\in Z^*$ such that $\sigma(z^*|A)<z^*(a+z)$ is satisfied for some $a\in A$. As $z^*(a+z)\leq \sigma(z^*|A)+z^*(z)$, this implies
$-z^*(z)< 0$ and $-\sigma(z^*|A)\neq-\infty$. But as $C\subseteq 0^+A$, $-\sigma(z^*|A)\neq-\infty$ implies $z^*\in\C$. Especially, $z$ is not an element of the right hand side of \eqref{eq:RecCone_finiteScal_2}.

On the other hand, assume $z\in 0^+A$, then $A$ is nonempty and $A+\cb{z}\subseteq A$, hence for all $z^*\in Z^*$ it holds
$\sigma(z^*|A+\cb{z})\leq \sigma(z^*|A)$, hence $\sigma(z^*|A)+z^*(z)\leq \sigma(z^*|A)$. This implies that either $-\sigma(z^*|A)=-\infty$ or $0\leq-z^*(z)$ is true for all $z^*\in Z^*$ and thus especially for $z^*\in\C$.

\pend

\begin{lemma}
Let $A\in\G^\triup$ be a nonempty set, then
\begin{equation*}
\cb{z^*\in \C\st -\sigma(z^*|A)\in\R}\subseteq (0^+A)^-\subseteq C^-.
\end{equation*}
\end{lemma}
\proof
Since $C\subseteq 0^+A$ is always satisfied, the last inclusion is trivial. Now take $z^*\in\C$ such that
$-\sigma(z^*|A)\in\R$ and $z\in 0^+A$, i.e. $A+z\subseteq A$. Then
\[
-\sigma(z^*|A)\leq -\sigma(z^*|A+z)=-\sigma(z^*|A)-z^*(z)
\]
implies $0\leq -z^*(z)$, in other words $z^*\in (0^+A)^-$.

\pend

%{\color{red}
For future reference, we collect some results on recession cones.

\begin{lemma}\label{lem:RecConeCalculus}
  Let $A,B\in\G^\triup$, $\mathcal A\subseteq \G^\triup$ and $s>0$ be given.
  \begin{enumerate}[(a)]
    \item
    It holds $0^+(sA)=0^+(A)$;
    \item
    If both $A$  and $B$ are nonempty, then $0^+(A+B)\subseteq 0^+(A)+0^+(B)$;
    \item
    If $A\idif B$ is nonempty, then $0^+(A)\subseteq 0^+(A\idif B)$;
  %  \item
  %  $\cl 0^+(A)\subseteq 0^+(\cl A) $;
  %  \item
  %  $\co 0^+(A)\subseteq 0^+(\co A)$;
    \item
    If $\bigcap\limits_{A\in\mathcal A}A$ is nonempty, then $ \bigcap\limits_{A\in \mathcal A}0^+(A)\subseteq 0^+(\bigcap\limits_{A\in \mathcal A}A)$;
    \item
    $\bigcap\limits_{A\in \mathcal A}0^+(A)\subseteq 0^+(\cl\co\bigcup\limits_{A\in \mathcal A}A) $.
    \end{enumerate}
\end{lemma}
\proof
  \begin{enumerate}[(a)]
    \item
    %It holds $0^+(tA)=0^+(A)$;
    $z\in 0^+(A)$ implies $a+\frac{t}{s}z\in A$ for all $a\in A$, hence $sa+tz\in sA$ is true for all $a\in A$;
    \item
    %If both $A$  and $B$ are nonempty, then $0^+(A+B)\subseteq =0^+(A)+0^+(B)$;
    Let $a\in A$ and $b\in B$ be given, $z\in 0^+(A)$, then $a+b+tz\in A+B$ is true for all $t>0$, hence $0^+(A)$ and likewise $0^+(B)$ is a subset of $0^+(A+B)$;
    \item
    %If $A\idif B$ is nonempty, then $0^+(A\idif B)\supseteq 0^+(A)$;
    Let $k\in A\idif B$, i.e. $B+k\subseteq A$ and $z\in 0^+(A)$, then $B+k+tz\subseteq A$ is true for all $t>0$, hence $0^+(A)\subseteq 0^+(A\idif B)$;
%    \item
%    %$0^+(\cl A)\supseteq \cl 0^+(A)$;
%    Let $\bar a\in \cl A$ be given, then for any $0$-neighborhood $U$ there is an $a\in (\bar a+U)\cap A$ and for any $z\in 0^+(A)$ it holds $a+tz\in (\bar a+tz)+U\cap A$, hence $\bar a+tz\in\cl A$ is true for all $t>0$, proving $0^+(A)\subseteq 0^+(\cl A)$. But as $0^+(\cl A)$ is closed, this already proves the statement;
%    \item
%    %$0^+(\co A)\supseteq \co 0^+(A)$;
%    Let $a,b\in A$ be given, $z\in 0^+(A)$ and $s\in\sqb{0,1}$. Then $sa+(1-s)b+stz+(1-s)tz\in A$ is true and as $0^+(\co A)$ is convex, this proves the statement;
    \item
    %If $\bigcap\limits_{A\in\mathcal A}A$ is nonempty, then $0^(\bigcap\limits_{A\in \mathcal A}A)\supseteq \bigcap\limits_{A\in \mathcal A}0^+(A)$;
    Let $a\in \bigcap\limits_{A\in\mathcal A}A$ be given, $z\in \bigcap\limits_{A\in \mathcal A}0^+(A)$, then for all $t>0$ it holds $a+tz\in \bigcap\limits_{A\in\mathcal A}A$, hence $z\in 0^+(\bigcap\limits_{A\in\mathcal A}A)$;
    \item
    %$0^+(\cl\co\bigcup\limits_{A\in \mathcal A}A)\supseteq \bigcap\limits_{A\in \mathcal A}0^+(A)$.
    Let $a\in \cl\co\bigcup\limits_{A\in \mathcal A}A$ be given.
    Then to any $U\in\mathcal U_Z$ there exist $a_1\in A_1$, $a_2\in A_2$ with $A_1, A_2\in\mathcal A$ and $s\in\sqb{0,1}$ such that $sa_1+(1-s)a_2\in a+U$ is true. Let $z\in\bigcap\limits_{A\in\mathcal 0^+A}$ be given, then especially $z\in 0^+A_1$ and $z\in 0^+ A_2$ is true and for any $t>0$ it holds
    \[
    sa_1+(1-s)a_2+tz\in \co\bigcup\limits_{A\in \mathcal A}A\cap(a+tz+U).
    \]
    This implies $a+tz\in\cl\co\bigcup\limits_{A\in\mathcal A}A$ is true for all $t>0$, hence by definition $z\in 0^+\cl\co\bigcup\limits_{A\in\mathcal A}A$.
    \end{enumerate}

\pend

%}

In the following proposition, we state some implications that are used in the main proofs.

\begin{proposition}\label{prop:A_nsubset_intB}%\label{lem:A_nsubset_intB_2}
Let $A, B\in\G^\triup$ be two sets, then
\begin{enumerate}[(a)]
\item
$A\nsubseteq \Int B$
 implies
\item
$\exists z^*\in W^*:\quad -\sigma(z^*|A)\leq-\sigma(z^*|B)\neq-\infty$
which in turn implies
\item
$\forall U\in\mathcal U_Z:\quad A\oplus U\nsubseteq B$.
\end{enumerate}
\end{proposition}
\proof
As $\Int C\neq \emptyset$, $B=\emptyset$ is equivalent to $\Int B=\emptyset$. In either case, $A\nsubseteq \Int B$ implies $A\neq\emptyset$, hence by a separation argument the inequality $ -\sigma(z^*|A)\leq-\sigma(z^*|B)\neq-\infty$ is satisfied for at least one $z^*\in W^*$.
%Otherwise, $\Int B\neq\emptyset$ is a convex set and by a separation argument there exists a $a\in A$ such that $-z^*(a)\leq -\sigma(z^*|B)$ is true for some $z^*\in Z^*\setminus\cb{0}$. But this implies $-\sigma(z^*|B)\neq-\infty$ and without loss $z^*\in W^*$.

For the second implication, consider that
for any $U\in\mathcal U_Z$ and any $z^*\in W^*$, $-\sigma(z^*|A+U)=\of{-\sigma(z^*|A)}\isum\of{-\sigma(z^*|U)}$ and by \cite[Remark 3.32]{HeydeSchrage11R} there exists a $\mu>0$ such that for all $w^*\in W^*$ it holds $-\sigma(w^*|U)\leq-\mu$.
Especially, if $-\sigma(z^*|A)\leq-\sigma(z^*|B)\neq-\infty$, then
\[
-\sigma(z^*|A+U)\leq-\sigma(z^*|A)-\mu < -\sigma(z^*|B),
\]
implying $A+U\nsubseteq B$.
\pend

The reverse implications do not hold in general, as the following example shows.

\begin{example}\label{ex:Inclusions}
\begin{enumerate}[(a)]
\item
Let $Z=\R^2$ and $C=B=\R^2_+$. Setting $A=\cb{(x,y)\in Z\st \frac{1}{x}\leq y,\, 0<x}$, then $A\subseteq \Int B$ but $-\infty=-\sigma((0,1)^T|A)<-\sigma((0,1)^T|B)$ is true.
\item
Let $Z=\R^2$ and $C=\cl\cone\cb{(0,1)^T}$. Set $A=\cb{(x,y)\in  Z\st x^2\leq y}$ and
\[
\forall n\in\N:\quad B_n=\cb{(x,y)\in Z\st \max\cb{2nx-n^2-\frac{1}{n},-2nx-n^2-\frac{1}{n}}\leq y },
\]
and $B=\bigcap\limits_{n\in\N}B_n$.
Then for all $U\in\mathcal U_Z$, $A+U\nsubseteq B$ is satisfied, but
 $-\sigma(z^*| B)=-\sigma(z^*|A)=-\infty$ is satisfied for $z^*\in\cone\cb{(-1,0)^T, (0,-1)^T}$ while for all other $z^*\in W^*$ it holds  $-\sigma(z^*| B)<-\sigma(z^*|A)$.
\end{enumerate}
\end{example}

In the reminder of this section we recall some properties and basic results about functions mapping into $\G(Z,C)$. First we need to point out that the conlinear space structure of the image space allows for intuitive definitions of properties of set--valued maps, namely straightforward extensions from scalar ones. Indeed a function $f \colon X \to \G^\triup$ is called convex when
\begin{equation*}
\label{EqConvFct} \forall x_1, x_2 \in X, \; \forall t \in \of{0,1} \colon
 f\of{tx_1+(1-t)x_2} \supseteq tf\of{x_1} \isum \of{1-t}f\of{x_2}.
\end{equation*}
Since $\supseteq$ stands for $\leq$, this definition mocks the scalar one. Moreover, $f$ is positively homogeneous (see e.g. \cite{Hamel09}) when
\[
\forall 0<t, \forall x \in X \colon f\of{tx} = tf\of{x},
\]
and it is called sublinear if it is positively homogeneous and convex.

Given a function $f:X\to \G^\triup$ we denote the (effective) domain by the set $\dom f=\cb{x\in X\st f(x)\neq \sup \G^\triup}$. Especially, if $\psi:S\subseteq X\to Z$ and  $f(x)=\psi^C(x)$ for all $x\in X$, then  $\dom f=S$, even though $f$ is defined on the whole set $X$. The image set of a subset $A\subseteq X$ through $f$ is denoted by $f\sqb{A}=\cb{f(x)\in \G^\triup\st x\in A}\subseteq \G^\triup$. A function $f:X\to \G^\triup$ is called proper, if $\dom f\neq\emptyset$ and $\inf \G^\triup\notin f\sqb{X}$.\\

 %In Section 3, we motivate the choice of $t\not=0$.\\

%We prove some more properties of set--valued functions $f \colon X \to \G^\triup$.
%For one thing, if $f(x)=\psi^C(x)$ for all $x\in X$, then  $\dom \psi^C=S$, even though $f$ is defined on the whole set $X$.

We introduce the family of extended real-valued functions $\vp_{f, z^*} \colon  X \to \R\cup\cb{\pm\infty}$ defined by
\[
\forall z^* \in \C:\quad \vp_{f,z^*}\of{x} = \inf\cb{-z^*\of{z} \mid z \in f\of{x}}
\]
as the family of scalarizations for $f$.
Some properties of $f$ are inherited by its scalarizations and vice versa. For instance,
$f$ is convex if and only if  $\vp_{f, z^*}$ is convex for each $z^* \in W^*$.
Moreover, by \eqref{eq:scal_representation_Set} the following representation is immediate
\[
\forall x \in X:\quad f\of{x} = \bigcap_{z^* \in W^*}\cb{z \in Z \st \vp_{f, z^*}\of{x} \leq -z^*\of{z}}.
\]

To state our main results, we need a notion of lower semicontinuity of set--valued functions. The following definition recalls some notions previously used in the literature, compare \cite{HeydeSchrage11R, Luc89, Zalinescu02}, to name but a few.

\begin{definition}\label{def:l.s.c.}
\begin{enumerate}[(a)]
\item
Let $\vp:X\to\OLR$ be a function, $x_0\in X$. Then $\vp$ is said to be lower semicontinuous (l.s.c.) at $x_0$, iff
\begin{equation*}
\forall r\in\R:\quad r<\vp_{f,z^*}(x_0)\;\Rightarrow\;
\exists U\in\mathcal U_X:\, \forall u\in U:\,  r<\vp_{f,z^*}(x_0+u).
\end{equation*}

\item
Let $f:X\to\G^\triup$ be a function, $M^*\subseteq\C$. Then $f$ is said $M^*$-- lower semicontiuous ($M^*$--l.s.c.) at $x_0$, iff
$\vp_{f,z^*}$ is l.s.c. at $x_0$ for all $z^*\in M^*$.
%
%\item
%Let $f:X\to\G^\triup$ be a function, $x_0\in X$ and $M^*\subseteq\C$. Then $f$ is said lower $M^*$--equicontinuous at $x_0\in\dom f$, iff
%\begin{align}\label{uniflsc}
%\forall \varepsilon >0 \;\exists U\in \mathcal U(x_0) \; \forall x\in U\; \forall z^*\in W^* : \quad  \varphi_{(f,z^*)}(x)>\varphi_{(f,z^*)}(x_0)\idif\varepsilon;
%\end{align}
%
%\item
%Let $\psi:S\subseteq X\to Z$ be a function, then $\psi$ is called $C$--continuous at $x_0\in S$, iff
%\[
%\forall V\in \mathcal U(\psi(x_0))\;\exists U\in\mathcal U(x_0)\;\forall x\in S\cap U:\quad \psi(x)\in V+C;
%\]

\item
Let $f:X\to\G^\triup$ be a function. If
\[
f(x)\supseteq \liminf\limits_{u\to 0}f(x+u)=\bigcap\limits_{U\in\mathcal U_X}\cl\co \bigcup\limits_{u\in U}f(x+u)
\]
is satisfied, then $f$ is lattice lower semicontinuous (lattice l.s.c.) at $x$.
%A function $f:X\to\G^\triup$ is lattice l.s.c. if and only if it is lattice l.s.c. everywhere.
\end{enumerate}
\end{definition}
In \cite{HeydeSchrage11R}, it has been proven that if $f$ is $\C$--l.s.c. at $x$, then it is also lattice l.s.c. at $x$.
Since we assume $\Int C\neq\emptyset$, $f$ is $\C$--l.s.c. at $x$ if and only if $f$ is $W^*$--l.s.c. at $x$.
One can show that if $f$ is convex, then $f$ is lattice l.s.c.
if and only if $\gr f=\cb{(x,z)\st z\in f(x)} \subseteq X \times Z$ is a closed set with respect to the product
topology, see \cite{HamelSchrage13PJO}.

Finally, we come back to weak efficiency. Obviously $x\in S$ is a weak solution to \eqref{VOP} if and only if one of the following equivalent assumptions is satisfied.
\begin{enumerate}[(a)]
\item $\forall u\in X:\quad \psi^C(x)\nsubseteq\Int\psi^C(x+u)$;
\item $\forall u\in X,\,\exists z^*\in W^*:\quad -\sup\cb{z^*(z)|z\in\psi^C(x)}\leq -\sup\cb{z^*(z)|z\in\psi^C(x+u)}\neq-\infty$;
\item $\forall u\in X,\,\forall U\in\mathcal U_Z:\quad \psi^C(x)+U\nsubseteq \psi^C(x+u)$.
\end{enumerate}

\begin{remark}
We note that, while
$$-\sup\cb{z^*(z)|z\in\psi^C(x)}=\left\{\begin{array}{ll}
		-z^*\psi(x)\in\R & \mbox{\rm if\;} x\in S\\
		+\infty & \mbox{\rm elsewhere,}
\end{array}\right.
$$
is true for $\psi:S\subseteq X\to Z$, considering a more general set--valued function $f:X\to\P(Z)$, it may happen that $f\left(x\right)=Z$ or $-\sup\cb{z^*(z)|z\in f(x)}=-\infty$ occurs.
\end{remark}

When considering any set--valued function $f:X\to\G^\triup$ and the related (weak) optimization problem
\begin{align}\label{P}\tag{P}
\min f(x),\; x\in X.
\end{align}
a point $x_0\in\dom f$ is called a weak minimizer of $f$ when

\begin{align}\tag{W-Min}\label{eq:WMin}
f(x)=Z\quad\vee\quad \forall x\in X\,\forall U\in\mathcal U_Z:\quad f(x_0)\oplus U\nsubseteq f(x).
\end{align}

This notion of solution can be related to others known in the literature. In \cite[Definition 2.4 (1)]{HernandezMarin05}, weak l--minimal elements of a set $\mathcal A\subseteq \P(Z)\setminus\cb{\emptyset}$ are those elements $A\in \mathcal A$, such that for all $B\in\mathcal A$, $A\subseteq B+\Int C$ implies $B\subseteq A+\Int C$.
If $\mathcal A\subseteq \G^\triup$, $\Int C\neq\emptyset$ implies $\Int B=B+\Int C$ for all $B\in\G^\triup$. Hence $A$ is a weak l--minimal element of $\mathcal A\subseteq\G^\triup$, if for all $B\in\mathcal A$, $A\subseteq \Int B$ implies $B\subseteq \Int A$. Thus, either $A=Z$, or $A$ is weak l--minimal in $\mathcal A\subseteq\G^\triup$, if and only if there exists no $B\in\mathcal A$ such that $A\subseteq \Int B$. Therefore for any $x_0\in\dom f$, $f(x_0)$ is a weak l--minimal element of $f\sqb{X}$ if and only if
\begin{align}\tag{W-l-Min}\label{eq:WMin_alternative}
f(x)=Z\quad\vee\quad \forall x\in X:\quad f(x_0)\nsubseteq \Int f(x).
\end{align}
Applying Proposition \ref{prop:A_nsubset_intB} it easily follows that \eqref{eq:WMin_alternative} implies
\begin{align}\tag{Sc-W-Min}\label{eq:ScWMin}
f(x_0)=Z\quad\vee\quad \forall x\in X\,\exists z^*\in W^*:\quad \vp_{f,z^*}(x_0)\leq \vp_{f,z^*}(x)\neq-\infty,
\end{align}
which in turn implies \eqref{eq:WMin}. While in general none of these implications can be reverted, we have some advantages when $f=\psi^C$ is the $\G^\triup$--valued extension of a vector--valued function $\psi:S\subseteq X\to Z$.

\begin{proposition}\label{rem:wEff_wMin}
Let $\psi:S\subseteq X\to Z$ be a vector--valued function.
For $f=\psi^C:X\to\G^\triup$ and $x_0\in S$, the properties \eqref{eq:WMin}, \eqref{eq:WMin_alternative} and  \eqref{eq:ScWMin} are equivalent and satisfied if and only if $x_0$ is a weakly efficient solution of \eqref{VOP}.
\end{proposition}

\proof
We only need to proof that \eqref{eq:WMin} implies \eqref{eq:WMin_alternative}.
When $f=\psi^C$ and $x_0$ is a weak minimizer,  then by definition for all $x\in X$ and all $U\in {\mathcal U_Z}$ it holds
$$
\psi\left(x_0\right)+C+U\not\subseteq \psi\left(x\right)+C,
$$
as in this case $f(x_0)=Z$ is not possible.
Assume for some $x\in S$  it holds
$$
\psi(x_0)+C\subseteq \Int\of{ \psi(x)+C}=\psi(x)+\Int C,
$$
or equivalently
$$
\of{\psi(x_0)-\psi(x)}+C\subseteq \Int C,
$$
implying $\of{\psi(x_0)-\psi(x)}\in \Int C$. But this is true, if and only if it exists an $U\in\mathcal U_Z$ such that
$$
\psi(x_0)+U\subseteq \psi(x)+\Int C,
$$
again implying
$$
\psi(x_0)+C+U\subseteq \psi(x)+C+\Int C\subseteq \psi(x)+C,
$$
a contradiction.
\pend

\begin{remark}\rm
For notational simplicity we set the restriction of a set--valued function $f:X\to\G^\triup$ to a segment with end points $x_0,x\in X$ as
$f_{x_0,x}:\R\to\G^\triup$, given by
\[
f_{x_0,x}(t)=\begin{cases}
	f(x_0+t(x-x_0)),\text{ if } t\in\sqb{0,1};\\
	\emptyset,\text{ elsewhere}
\end{cases}
\]
and the restriction of a scalar--valued function $\vp:X\to\OLR$ to the same segment by
\[
\vp_{x_0,x}(t)=\begin{cases}
	\vp(x_0+t(x-x_0)),\text{ if } t\in\sqb{0,1};\\
	+\infty,\text{ elsewhere.}
\end{cases}
\]
Setting $x_t=x_0+t(x-x_0)$ for all $t\in\R$,
the scalarization of the restricted function $f_{x_0,x}$ is equal to the restriction of the scalarization of $f$ for all $z^*\in\C$.

If $f$ is convex, $x_0,x_t\in \dom f$ for some $t\in\of{0,1}$, then $\of{\vp_{f,z^*}}_{x_0,x}$ is lower semicontinuous on $\of{0,t}$ for all $z^*\in \C$, hence $f_{x_0,x}$ is $\C$--l.s.c. on $\of{0,t}$.

Notice that in general, if $f$ is $\C$--l.s.c. in $x_0$, then $f_{x_0,x}$ is $\C$--l.s.c. in $0$ for all $x\in X$, while the implication is not revertible.
\end{remark}

%%%%%%%%%%%%%%%%%%%%%%%%%%%%%%%%%%%%%%%%%%%%%%%%%%%

\section{Dini Directional Derivatives}\label{sec:DirDer}

As we anticipated in Section \ref{sec:Setting}, $\inf$--residuated and order complete  structure allows for an immediate extension of the definitions of both the difference quotient and upper and lower limits. Thus we have the basic ingredients to define the notion of upper and lower Dini directional derivatives.

\begin{definition}\label{def:DirDer}
Let  $f:X\to \G^\triup$ and $x,u\in X$. The {\em upper} and {\em lower Dini directional derivative} of $f$ at $x$ in direction $u$ are given by
\begin{align*}
f^\uparrow(x,u)&=\limsup\limits_{t\downarrow 0}\frac{1}{t}\of{f(x+tu)\idif f(x)}
=\inf\limits_{0<s}\sup\limits_{0<t\leq s}\frac{1}{t}\of{f(x+tu)\idif f(x)};\\
f^\downarrow(x,u)&=\liminf\limits_{t\downarrow 0}\frac{1}{t}\of{f(x+tu)\idif f(x)}
=\sup\limits_{0<s}\inf\limits_{0<t\leq s}\frac{1}{t}\of{f(x+tu)\idif f(x)}.
\end{align*}
If both derivatives coincide, then $f'(x,u)=f^\uparrow(x,u)=f^\downarrow(x,u)$ is the {\em Dini directional derivative} of $f$ at $x$ in direction $u$.
\end{definition}
It is easy to see that $f^\downarrow(x,u)\supseteq f^\uparrow(x,u)$ is always satisfied, hence $f'(x,u)$ exists if and only if $f^\uparrow(x,u)\supseteq f^\downarrow(x,u)$.

The previous definition does not require $f$ to be proper or $x\in\dom f$. Clearly, if $f(x)=\sup \G^\triup=\emptyset$, then $f'(x,u)=\inf \G^\triup=Z$ is satisfied for all $u\in X$.

%For notational simplicity, we agree that in the sequel we refer simultaneously to upper and lower Dini derivatives by $f^\updownarrow(x,u)$, even if the two values can be different.

If $0<s$ is given, then $f^\uparrow(x,su)=s f^\uparrow(x,u)$ and $f^\downarrow(x,su)=s f^\downarrow(x,u)$, that is, both derivatives are positively homogeneous in the second component.
	We also remark that, when $u=0$, $\frac{1}{t}\of{f(x+t0)\idif f(x)}\supseteq \frac{1}{t}\cdot C$, so both derivatives in direction $0$ may not be equal to 
 $C$, the neutral element in $\G^\triup$. This explains the choice not to include the assumption $f(0)=C$ in the the definition of positive homogeneity in the previous section.

\begin{remark}\rm
When $Z=\R$, Definition \ref{def:DirDer} provides an extension to the classical notion of Dini derivatives for scalar functions (see \cite{giorgi1992dini} and the references therein), without requiring neither $x\in\dom f$ nor $f$ to be proper.
However, since a vector space needs not to be order complete, the same definition may not be applied to vector--valued functions $\psi:S\subseteq X\to Z$.
For this reason, different Dini derivatives for vector--valued functions have been defined, compare e.g. \cite{AnsariLee2010, Ginchev2007}.
\end{remark}

\begin{example}\label{ex:scalar_DirDer}
Let $\vp:X\to\OLR$ be an extended scalar function. If $\vp(x+tu)\in\R$ is satisfied for all $t\in\sqb{0,t_0}$ for a given $0<t_0$, then the difference quotient is real
$$
\displaystyle\frac{1}{t}\left(\vp\left(x+tu\right)-\vp\left(x\right)\right)\in\R.
$$
Hence in this case the Dini derivatives coincide with the standard definition in the literature, compare \cite{giorgi1992dini}.\\
If $x\notin\dom \vp$, then $\vp(x+tu)\idif \vp(x)=-\infty$ for all $t>0$, so $\vp'(x,u)=-\infty$.
On the other hand, if $\vp(x)=-\infty$, then $\vp(x+tu)\idif \vp(x)=-\infty$, whenever $\vp(x+tu)=-\infty$ and $\vp(x+tu)\idif \vp(x)=+\infty$, else. The value of the derivatives in this case depends on the behaviour of $\vp$ in a proximity of $x$.
\end{example}

\begin{proposition}\label{prop:DirDer_of_conv}
Let  $f:X\to \G^\triup$ be convex, then the Dini derivative exists for all $u\in X$ and it holds
\[
f'(x,u)=\inf\limits_{t\downarrow 0}\frac{1}{t}\of{f(x+tu)\idif f(x)}.
\]
Moreover, $f':X\times X\to \G^\triup$ is sublinear in its second component.
\end{proposition}
\proof
Let $0<s$ be given, then for all $0<t\leq s$, there exists a $0<h\leq 1$ such that $hs+(1-h)0=t$ and by convexity of $f$, $f(x+tu)=f(h(x+su)+(1-h)x)\supseteq hf(x+su)\isum (1-h)f(x)$. By assumption,  $hf(x)\isum (1-h)f(x)\supseteq f(x)$ is satisfied for all $h\in\sqb{0,1}$.
Applying Lemma \ref{lem:calc_conlin} we can prove
\begin{align*}
\frac{1}{t}\of{f(x+tu)\idif f(x)}
&\supseteq \frac{1}{hs}\of{\of{hf(x+su)\isum (1-h)f(x)}\idif \of{hf(x)\isum (1-h)f(x)}}\\
&\supseteq \frac{1}{hs}\of{\of{hf(x+su)\idif hf(x)}\isum \of{(1-h)f(x)\idif (1-h)f(x)}}\\
&=   \frac{1}{hs}\of{h\of{f(x+su)\idif f(x)}\isum (1-h)\of{f(x)\idif f(x)}}\\
&\supseteq \frac{1}{hs}\of{h\of{f(x+su)\idif f(x)}\isum \theta}\\
&= \frac{1}{s}\of{f(x+su)\idif f(x)}.
\end{align*}
Hence especially
\[
f^\uparrow(x,u)\supseteq \inf\limits_{0<s}\frac{1}{s}\of{f(x+su)\idif f(x)}\supseteq f^\downarrow(x,u)
\]
is proven.

Finally, let $s\in\of{0,1}$ and $0<t\leq r$ be given, $u_1, u_2\in X$. Then
\begin{align*}
f'(x,su_1+(1-s)u_2)
&\supseteq \frac{1}{t}\of{f(x+t(su_1+(1-s)u_2))\idif f(x)}\\
&\supseteq \frac{1}{t}\of{s\of{f(x+tu_1)\idif f(x)}\isum (1-s)\of{f(x+tu_2)\idif f(x)}}\\
&\supseteq s\frac{1}{t}\of{f(x+tu_1)\idif f(x)}\isum (1-s)\frac{1}{r}\of{f(x+ru_2)\idif f(x)}
\end{align*}
But as this holds for all $0<t\leq r$,
\begin{align*}
f'(x,su_1+(1-s)u_2)\supseteq sf'(x,u_1)\isum (1-s)\frac{1}{r}\of{f(x+ru_2)\idif f(x)}
\end{align*}
and ultimately
\begin{align*}
f'(x,su_1+(1-s)u_2)\supseteq sf'(x,u_1)\isum (1-s)f'(x,u_2)
\end{align*}
are true.
As $f'(x,\cdot):X\to \G^\triup$ is convex and positively homogeneous, it is sublinear.

\pend

In the proof of Proposition \ref{prop:DirDer_of_conv}, it is shown that under the given assumptions the difference quotient is decreasing.

%If $Y=\G^\triup$, the assumption $tf(x)\isum (1-t)f(x)\lel f(x)$ is always satisfied, because of the ordering relation given by $\supseteq$. Thus, if $f:X\to\G^\triup$ is convex, then $tf(x)\isum (1-t)f(x)= f(x)$ is satisfied for all $x\in X$ and all $t\in\of{0,1}$.

\begin{lemma}\label{prop:DirDer_of_conv_SV}
Let $f:X\to\G^\triup$ be a convex function, $x,u\in X$. Then
\[
f'(x,u)=\cl\bigcup\limits_{0<t\leq s}\frac{1}{t}\of{f(x+tu)\idif f(x)}
\]
is true for all $s>0$. Moreover it holds
\[
\Int f'(x,u)=\bigcup\limits_{0<t\leq s}\Int \frac{1}{t}\of{f(x+tu)\idif f(x)}.
\]
\end{lemma}
\proof
Proposition \ref{prop:DirDer_of_conv} proves $f'(x,u)=\inf\limits_{0<t\leq s}\frac{1}{t}\of{f(x+tu)\idif f(x)}$.
Moreover, since the difference quotient is decreasing as $t$ converges to $0$, $\bigcup\limits_{0<t\leq s}\frac{1}{t}\of{f(x+tu)\idif f(x)}$ is convex for all $0<s$, the first statement is true.

As for the second statement, let $z\in \Int f'(x,u)$ be given. Especially
there exists  $\bar z\in\Int C$ and $U\in\mathcal U_Z$ such that $z-\bar z\in\Int f'(x_0,x)$ and $\bar z+U\subseteq \Int C$.
 Therefore, there exists $0<t$ such that $\of{z-\bar z}\in \frac{1}{t}\of{f(x+tu)\idif f(x)}$ and
\[
z-\bar z +\bar z+U\subseteq \frac{1}{t}\of{f(x+tu)\idif f(x)}+\Int C\subseteq \frac{1}{t}\of{f(x+tu)\idif f(x)},
\]
implying $z+U\subseteq \frac{1}{t}\of{f(x+tu)\idif f(x)}$, or equivalently $z\in \Int \frac{1}{t}\of{f(x+tu)\idif f(x)}$.
\pend

\begin{proposition}\label{prop:Rec_at_x_at_f'(x,u)}
Let $f:X\to\G^\triup$ be given, $x,u\in X$.  If $f^\uparrow(x,u)\neq\emptyset$ ($f^\downarrow(x,u)\neq\emptyset$), then
$0^+f(x)\subseteq 0^+f^\uparrow(x,u)$ ($0^+f(x)\subseteq 0^+f^\downarrow(x,u)$)
is true. If additionally  $x\in \dom f$, then $f'(x,0)=0^+f(x)$.
\end{proposition}
\proof
First we consider the case $x\in \dom f$ and $u=0$, then
\[
f'(x,0)%=\liminf\limits_{t\downarrow 0}\frac{1}{t}\of{f(x)\idif f(x)}
=f(x)\idif f(x)=0^+f(x).
\]
For arbitrary $u\in X$, let $f^\uparrow(x,u)\neq \emptyset$ be satisfied.
By Lemma \ref{lem:RecConeCalculus}, it holds
\[
0^+(f^\uparrow(x,u))\supseteq \cl\co\bigcup\limits_{s>0}\bigcap\limits_{t\in\of{0,s}}0^+(f(x+tu)\idif f(x))\supseteq\cl\co 0^+(f(x))
\]
and as $f(x)$ is itself closed and convex, this implies $0^+(f^\uparrow(x,u))\supseteq 0^+(f(x))$.
By the same argument,
\[
0^+(f^\downarrow(x,u))\supseteq \bigcap\limits_{s>0}\cl\co\bigcup\limits_{t\in\of{0,s}}0^+(f(x+tu)\idif f(x))\supseteq 0^+(f(x))
\]
holds true.

%This is equivalent to stating that it exists an $s>0$ such that for all $0<t\leq s$ it holds $f(x)+\cb{tz}\subseteq f(x+tu)$. But, as $f(x)$ can be re-written as $f(x)+t\cdot0^+f(x)$ for all $0<t$, the previous inclusion is equivalent to
%\[
%\forall 0<t\leq s:\quad \cb{z}+0^+f(x)\subseteq \frac{1}{t}\of{f(x+tu)\idif f(x)},
%\]
%hence $0^+f(x)\subseteq 0^+f^\uparrow(x,u)$.
%Again, the argument for the lower derivative goes along the same lines.
\pend

Especially, if $f$ is convex, $x\in \dom f$, then $f'(x,\cdot):X\to \G^\triup$ is a sublinear function with $f'(x,0)=0^+f(x)$, the neutral element in a subspace of the image space. However, $0^+f(x)\supseteq C$ and in general, the inequality will be strict.

We are also interested in comparing the derivative of a given function with the set of the derivatives of its scalarization. The following inequalities holds true.

\begin{proposition}\label{prop:vp'fLeqf'}
Let $f:X\to\G^\triup$ be given, $x,u\in X$ and $z^*\in W^*$. Then
\begin{align*}
f^\uparrow(x,u)&\subseteq \bigcap\limits_{z^*\in W^*}\cb{z\in Z\st \vp_{f,z^*}^\uparrow(x,u)\leq -z^*(z)};\\
\vp_{f,z^*}^\uparrow(x,u)	&\leq -\sigma(z^*|f^\uparrow(x,u))
\end{align*}
and 
\begin{align*}
f^\downarrow(x,u)&\subseteq \bigcap\limits_{z^*\in W^*}\cb{z\in Z\st \vp_{f,z^*}^\downarrow(x,u)\leq -z^*(z)};\\
\vp_{f,z^*}^\downarrow(x,u)	&\leq -\sigma(z^*|f^\downarrow(x,u)).
\end{align*}

\end{proposition}
 \proof
Combining the scalarization formula \eqref {eq:scal_representation_Set} with Lemmas \ref{prop:scal_of_infimum} to \ref{lem:scal_of_difference}, it holds
\begin{align*}
\vp_{f,z^*}^\uparrow(x,u)
&=\inf\limits_{0<s}\sup\limits_{0<t\leq s}\frac{1}{t}\of{\vp_{f,z^*}(x+tu)\idif \vp_{f,z^*}(x)}\\
&\leq\inf\limits_{0<s}\sup\limits_{0<t\leq s}-\sigma(z^*|\frac{1}{t}\of{f(x+tu)\idif f(x)})\\
&\leq\inf\limits_{0<s}-\sigma(z^*|\bigcap\limits_{0<t\leq s}\frac{1}{t}\of{f(x+tu)\idif f(x)})\\
&=-\sigma(z^*|\cl\co\bigcup\limits_{0<s}\bigcap\limits_{0<t\leq s}\frac{1}{t}\of{f(x+tu)\idif f(x)})\\
&=-\sigma(z^*|f^\uparrow(x,u)).
\end{align*}
This immediately implies
\begin{align*}
f^\uparrow(x,u)\subseteq \bigcap\limits_{z^*\in W^*}\cb{z\in Z\st \vp_{f,z^*}^\uparrow(x,u)\leq -z^*(z)}.
\end{align*}
%and likewise
%\begin{align*}
%f^\uparrow(x,u)&=\cl\co\bigcup\limits_{0<s}\bigcap\limits_{0<t\leq s}\frac{1}{t}\bigcap\limits_{z^*\in W^*}
%\cb{z\in Z\st \vp_{f,z^*}(x+tu)\idif \vp_{f,z^*}(x)\leq -z^*(z)}\\
%&\subseteq\bigcap\limits_{z^*\in W^*}\cl\co\bigcup\limits_{0<s}\bigcap\limits_{0<t\leq s}\frac{1}{t}
%\cb{z\in Z\st \vp_{f,z^*}(x+tu)\idif \vp_{f,z^*}(x)\leq -z^*(z)}\\
%&\subseteq\bigcap\limits_{z^*\in W^*}\cl\co\bigcup\limits_{0<s}
%\cb{z\in Z\st \sup\limits_{0<t\leq s}\frac{1}{t}\of{\vp_{f,z^*}(x+tu)\idif \vp_{f,z^*}(x)}\leq -z^*(z)}\\
%&=\bigcap\limits_{z^*\in W^*}
%\cb{z\in Z\st \inf\limits_{0<s}\sup\limits_{0<t\leq s}\frac{1}{t}\of{\vp_{f,z^*}(x+tu)\idif \vp_{f,z^*}(x)}\leq -z^*(z)}
%\end{align*}
The same chain of arguments proves both inequalities for the lower derivative as well.
\pend

In general, neither of the inequalities in  Proposition \ref{prop:vp'fLeqf'} is satisfied with equality, as the following counterexample shows. That is, the operations of taking the derivative and taking the scalarization of a function do not commute.

\begin{example}
Let $f:\R\to\G^\triup(\R,\cb{0})$ be defined as $f(x)=\sqb{-\sqrt{1-x^2}, \sqrt{1-x^2}}$, whenever $x\in\sqb{-1,1}$ and $f(x)=\emptyset$, else.
Then $f'(0)+z\nsubseteq f(t)$ for any $t\neq 0$, so $f'(0,u)=\emptyset$.
On the other hand, $\vp_{f,s}(x)=-|s|\cdot \sqrt{1-x^2}$ for all $s\neq 0$ and thus $\vp'_{f,s}(x,u)=-|s|\cdot \frac{x}{\sqrt{1-x^2}}\cdot u$ for all $x\in\of{-1,1}$, especially $\vp'_{f,s}(0,u)=0$ for all $s\neq 0$.
Hence,
\begin{align*}
\emptyset=f'(0,u)\subsetneq \bigcap\limits_{z^*\in(\cb{0})^-\setminus\cb{0}}f'_{z^*}(0,u)=\cb{0}
\end{align*}
\end{example}

\begin{example}
Let $\psi:S\subseteq X\to Z$ be a $C$--convex function with set--valued extension $f=\psi^C:X\to\G^\triup$, then for all $x,x+u\in S$ and all $t\in \of{0,1}$ it holds
\[
\frac{1}{t}\of{f(x+tu)\idif f(x)}=\frac{1}{t}\of{\psi(x+tu)-\psi(x)}+C.
\]
If $x\notin S$, then $f'(x,u)=Z$ while if $x\in S$ and $x+tu\notin S$ is satisfied for all $0<t$, then $f'(x,u)=\emptyset$.
Thus especially for $x\in S$,
\[
f'(x,u)=\cl\bigcup\limits_{\substack{0<t,\\ x+tu\in S}}\of{\frac{1}{t}\of{\psi(x+tu)-\psi(x)}+C}.
\]
is satisfied. However the infimum of $\frac{1}{t}\of{\psi(x+tu)-\psi(x)}$ needs not exist, even if $Z$ is lattice ordered.
\end{example}

Proposition \ref{prop:vp'fLeqf'} and the previous examples motivate to consider as a special case when equality is satisfied in either of the inequalities stated in Proposition \ref{prop:vp'fLeqf'}.
In the sequel we refer to
\begin{align}\label{eq:reg_sharp}\tag{SR}
\forall z^*\in B:\quad \vp_{f,z^*}^\downarrow(x,u)= -\sigma(z^*|f^\downarrow(x,u))
\end{align}
as strong regularity assumption, in contrast to the weak regularity assumption
\begin{align}\label{eq:reg_weak}\tag{WR}
f^\downarrow(x,u)&= \bigcap\limits_{z^*\in B}\cb{z\in Z\st \vp_{f,z^*}^\downarrow(x,u)\leq -z^*(z)}.
\end{align}

The following proposition states that if $f=\psi^C$, then it satisfies \eqref{eq:reg_weak}. Additionally assuming convexity allows to prove \eqref{eq:reg_sharp}.
It is left as an open question to identify necessary and sufficient conditions for either regularity assumption to be satisfied by a set--valued function $f:X\to\G^\triup$.

\begin{proposition}\label{prop:vp'fLeqf'2}
Let $\psi:S\subseteq X\to Z$ be given, $x,u\in X$ and $f=\psi^C:X\to\G^\triup$ its set--valued extension. Then property \eqref{eq:reg_weak} is satisfied for the lower derivative of $f$,
\begin{align*}
f^\downarrow(x,u)&= \bigcap\limits_{z^*\in B}\cb{z\in Z\st \vp_{f,z^*}^\downarrow(x,u)\leq -z^*(z)}.
\end{align*}
If additionally $\psi$ is $C$--convex, i.e. for all $x_1,x_2\in X$ and all $t\in\sqb{0,1}$ it holds
\[
t\psi(x_1)+(1-t)\psi(x_2)\subseteq \psi(tx_1+(1-t)x_2)+C
\]
 then property \eqref{eq:reg_sharp} is true.
\end{proposition}
 \proof
Recall that for all $z^*\in W^*$ it holds $\vp_{f,z^*}(x)=-z^*\psi(x)$ for all $x\in S$ and $\vp_{f,z^*}(x)=+\infty$, elsewhere. Hence especially
\[
-\sigma(z^*|\frac{1}{t}\of{f(x+tu)\idif f(x)})=\frac{1}{t}\of{\vp_{f,z^*}(x+tu)\idif \vp_{f,z^*}(x)}
\]
is satisfied for all $0<t$, in contrast to the inequality in the case of general set--valued functions.
Applying Lemma \ref{prop:scal_of_infimum}, then
\[
-\sigma(z^*|\inf\limits_{0<t}\frac{1}{t}\of{f(x+tu)\idif f(x)})=\inf\limits_{0<t}\frac{1}{t}\of{\vp_{f,z^*}(x+tu)\idif \vp_{f,z^*}(x)},
\]
proving the equality in the convex case. Also, additionally applying Lemma \ref{prop:scal_of_supremum}, it holds
\begin{align*}
f^\downarrow(x,u)&=\bigcap\limits_{0<s}\bigcap\limits_{z^*\in W^*}
\cb{z\in Z\st \inf\limits_{0<t\leq s}\frac{1}{t}\of{\vp_{f,z^*}(x+tu)\idif \vp_{f,z^*}(x)}\leq -z^*}\\
&=\bigcap\limits_{z^*\in W^*}
\cb{z\in Z\st \sup\limits_{0<s}\inf\limits_{0<t\leq s}\frac{1}{t}\of{\vp_{f,z^*}(x+tu)\idif \vp_{f,z^*}(x)}\leq -z^*}.
\end{align*}

\pend

\begin{example}
\begin{enumerate}[(a)]
\item
Let $Z=\R^3$ be ordered by the ordering cone $C$,  the closed conical hull of $\co\cb{(-1,1,1)^T,(-1,1,-1)^T,(1,1,-1)^T,(1,1,1)^T}$. Let $\psi:S\subseteq X\to Z$ be given with $\psi(0)=(0,0,0)^T$ and
\[
\psi(t)=\begin{cases}
(-t,0,0)^T,\text{ if } \exists n\in \N:\quad \frac{1}{2n}\leq t<\frac{1}{2n-1};\\
(t,0,0)^T,\text{ if } \exists n\in \N:\quad\frac{1}{2n+1}\leq t<\frac{1}{2n}.
\end{cases}
\]
Then $(\psi^C)^\uparrow(0,1)=\co\cb{(0,1,1)^T,(0,1,-1)^T}\oplus C$. For $z^*=(0,-1,0)\in \C$, it holds
$\vp'_{(\psi^C),z^*}(0,1)=0<-\sigma(z^*|(\psi^C)^\uparrow(0,1))=1$.

\item
Let $Z=\R^2$ be ordered by the natural ordering cone $C=\R^2_+$ and let a vector function $\psi:S\subseteq X\to Z$ be given such that $\psi(x)=(0,0)^T$ and
$\psi(t)=(1,0)^T$
is satisfied for all $t>0$. Then $(\psi^C)'(0,1)=\emptyset$, hence $-\sigma(z^*|(\psi^C)'(0,1))=+\infty$, while for $z^*=(0,-1)^T\in \C$ it holds $\vp'_{(\psi^C),z^*}(0,1)=0$.
\end{enumerate}
\end{example}

In Proposition \ref{prop:vp'fLeqf'2}, we basically apply set--valued arguments to obtain a definition of Dini derivatives for vector--valued functions. In \cite{AnsariLee2010, Ginchev2007}, similar derivatives are introduced using vector--valued arguments.\\
Although a careful comparison among the different types of derivatives is beyond the limits of the paper, we conclude this section with a sneak view of some results that easily hold. First, we stress once more that Definition \ref{def:DirDer} allows to  introduce a Dini type derivative without defining infinite elements in a vector space.

To compare our approach to that in  \cite{AnsariLee2010}, let $C$ be a polyhedral cone, $M^*\subseteq W^*$ a finite set such that $\co M^*= W^*$. If $\psi:S\subseteq X\to Z$ is a $C$--convex function, $f(x)=\psi^C(x)$ for all $x\in X$, then $\bar z\in\bigcap\limits_{m^*\in M^*}\cb{z\in Z\st \vp'_{f,m^*}(x,u)}$ implies that for all $t>0$ there exists $\varepsilon_t>0$ such that
\begin{align*}
\forall m^*\in M^*:\quad \frac{1}{t}\of{ \vp_{f,m^*}(x+tu)\idif \vp_{f,m^*}(x)}\leq -m^*(\bar z)+\varepsilon_t.
\end{align*}
As any $z^*\in W^*$ can be represented as a convex combination of elements of $M^*$, and $\vp_{f,z^*}(x)=-z^*\psi(x)$ for all $x\in S=\dom \vp_{f,z^*}$, this implies
\begin{align*}
\forall z^*\in W^*:\quad \frac{1}{t}\of{ \vp_{f,z^*}(x+tu)\idif \vp_{f,z^*}(x)}\leq -z^*(\bar z)+\varepsilon_t,
\end{align*}
hence $\bar z\in\bigcap\limits_{z^*\in W^*}\cb{z\in Z\st \vp'_{f,z^*}(x,u)\leq -z^*(z) }$. Therefore in this case,
\[
(\psi^C)'(x,u)=\bigcap\limits_{m^*\in M^*}\cb{z\in Z\st \vp'_{f,m^*}(x,u)}
\]
is satisfied, as $(\psi^C)'(x,u)\subseteq\bigcap\limits_{m^*\in M^*}\cb{z\in Z\st \vp'_{f,m^*}(x,u)}$ is always true.
Especially, if $Z=\R^n$ is ordered by the Pareto ordering cone, then the derivative of the set--valued extension of a $C$--convex function $\psi:S\subseteq X\to Z$ is characterized by the derivatives of the finite number of scalarizations with respect to the negative unit vectors in $Z^*$.
This approach has been chosen in \cite{AnsariLee2010}, where the upper and lower Dini derivative of a function $\psi:S\subseteq X\to\R^n$ is defined through the vector $(\vp^\downarrow_{\psi^C,-e_1^*}(x,u),...,\vp^\downarrow_{\psi^C,-e_n^*}(x,u))^T\in\OLR^n$, $e_i^*$ denoting the i-th unit vectors in $\R^n$.

In \cite{Ginchev2007}, a set--valued Dini derivative for vector--valued functions $\psi:S\subseteq X\to Z$ has been defined, using the Painelev\'e Kuratowski limit of the difference quotient. The original image space is extended by infinite elements $z_\infty=\lim\limits_{t\to\infty}tz$ for all $z\in Z\setminus\cb{0}$.
Roughly speaking, $z_\infty$ is an element of $\psi'(x,u)$, if for any $U\in\mathcal U_Z$ and any $s>0$, for any $t_0>0$  there exists a $t\in\of{0,t_0}$ such that
$\frac{1}{t}\of{\psi(x+u)\idif \psi(x)}\in sz+\cone(\cb{z}+U)$ and $z\in \psi'(x,u)$, if $z$ is a cluster point of the net of difference quotients.
It can be proven that if $z\in\psi'(x,u)$, then $z\in (\psi^C)^\downarrow(x,u)$, while the situation is somewhat more complicated for infinite elements.
If $z_\infty\in\psi'(x,u)$ and $z\in-\Int C$, then $(\psi^C)^\downarrow(x,u)=Z$.
With Stampacchia type variational inequalities in mind, the following chain of implications can be proven.

If $\psi'(x,u)\cap(-C\cup\cb{z_\infty\st z\in -C\setminus\cb{0}})=\emptyset$, then $0\notin (\psi^C)^\downarrow(x,u)$, which in turn implies $\psi'(x,u)\cap(-\Int C\cup\cb{z_\infty\st z\in \Int C})=\emptyset$.
%
%However, consider the function $\psi:\R\to\R^2$
%$$
%\psi(t)=
%\left\{\begin{array}{ll}
%(-1,-\sqrt t) & t>0\\
%(0,0) & t=0\\
%+\infty & \mbox{\rm elsewhere}
%\end{array}\right.
%$$
%and the image space $Z=\R^2$ being ordered by $\R^2_+$. Then the differential quotient of $\psi$ is $\frac{1}{t}\of{\psi(0+t)-\psi(0)}=(-\frac{1}{t},-\frac{1}{\sqrt t})$, so $(\psi^C)'(0,1)=Z$, while $\psi'(0,1)=\cb{(-1,0)_\infty}$.

%It turns out, by direct calculations, that also the function
%$$
%\Psi(t)=
%\left\{\begin{array}{ll}
%(-1,\sqrt t) & t>0\\
%(0,0) & t=0\\
%+\infty& \mbox{\rm elsewhere}
%\end{array}\right.
%$$
%then $(\Psi^C)'(0,1)=\emptyset$, while $\Psi'(0,1)=\cb{(-1,0)_\infty}$.
%

%%%%%%%%%%%%%%%%%%%%%%%%%%%%%%%%%%%%%%%%%%%%%%%%%%%%

\section{Main Results} \label{sec:Main_Results}

To characterize weak minimizers of \eqref{VOP} as solutions to (weak) variational inequalities of Stampacchia or Minty type, we first provide extensions of such inequalities for a general, convex, set--valued function $f:X\to\G^\triup$ and study their relations with solutions of \eqref{P}. \\

We begin by considering the following variational inequality of Stampacchia type.

\begin{definition}\label{def:SVI}
Let $f:X\to\G^\triup$ be a convex function and $f':X\times X\to\G^\triup$ its Dini derivative.
Then $x_0$ is  a solution to the {\em weak Stampacchia variational inequality}, iff
\begin{align}\tag{W-SVI}\label{eq:WSVI}
f(x_0)=Z\quad\vee\quad  \forall x\in X:\quad 0\notin\Int f'(x_0,x-x_0).
\end{align}
\end{definition}

\begin{remark}\label{rem:SVI}\rm
An element $x_0\in\dom f$ solves \eqref{eq:WSVI} if and only if
\begin{equation}\label{star}
f(x_0)=Z\quad\vee\quad  \forall x\in X:\quad 0^+f(x_0)\nsubseteq \Int f'(x_0,x-x_0)
\end{equation}
%or
%\[
%f(x_0)=Z\quad\vee\quad  \forall x\in X\,\forall U\in\mathcal U_Z:\quad U\nsubseteq f'(x_0,x-x_0)
%\]
is satisfied.

Indeed,
%$0\notin\Int f'(x_0,x-x_0)$ is equivalent to $f'(x_0,x-x_0)=\emptyset$ or $0^+\of{f'(x_0,x-x_0)}\nsubseteq \Int f'(x_0,x-x_0)$.
%Obviously,
$0\notin\Int f'(x_0,x-x_0)$ implies $0^+f(x_0)\nsubseteq \Int f'(x_0,x-x_0)$, as $x_0\in \dom f$ and hence $0\in 0^+f(x_0)$ is satisfied.
On the other hand $0\in\Int f'(x_0,x-x_0)$ implies
$0^+f(x_0)\subseteq 0^+f'(x_0,x-x_0)$ (compare Proposition \ref{prop:Rec_at_x_at_f'(x,u)}) and thus $0^+f(x_0)\subseteq \Int f'(x_0,x-x_0)$.
\end{remark}

According to the ordering relation introduced in $\G^{\triup}$, \eqref{star} can be easily read as an inequality in the conlinear space that perfectly matches the form of scalar variational inequalities.\\

Applying scalarization, we can prove relations between the set--valued inequality \eqref{eq:WSVI} the family of variational inequalities corresponding to the scalarizations of the given set--valued function.

\begin{lemma}\label{lem:SVI}
If $x_0\in\dom f$ satisfies
\begin{align}\tag{Sc-W-SVI}\label{eq:ScWSVI}
f(x_0)=Z\quad\vee\quad  \forall x\in X\,\exists z^*\in W^*:\; 0\leq\vp'_{f,z^*}(x_0,x-x_0)
\end{align}
 then it solves \eqref{eq:WSVI}.
If additionally the regularity assumption \eqref{eq:reg_sharp} is satisfied, the reverse implication is true, too.
\end{lemma}
\proof
By a separation argument, $0\notin\Int f'(x_0,x-x_0)$ is satisfied, if and only if there exists a $z^*\in W^*$ such that $0\leq-\sigma(z^*|f'(x_0,x-x_0))$. But as by Proposition \ref{prop:vp'fLeqf'} the inequality $\vp'_{f,z^*}(x_0,x-x_0)\leq -\sigma(z^*|f'(x_0,x-x_0))$ is always satisfied, the first implication is proven.
On the other hand if \eqref{eq:reg_sharp} is satisfied, then $\vp'_{f,z^*}(x_0,x-x_0)= -\sigma(z^*|f'(x_0,x-x_0))$ is true for all $z^*\in W^*$ and thus the reverse implication holds true.
\pend

Under convexity assumptions,  \eqref{eq:WSVI} is a necessary and sufficient condition for \eqref{eq:WMin} to hold.

\begin{theorem}\label{thm:Stamp_Min1}
Let $f:X\to G^\triup$ be a convex function, $x_0\in \dom f$.
Then $x_0$ is a weak minimizer of $f$ if and only if it solves the Stampacchia variational inequality \eqref{eq:WSVI}.
\end{theorem}
\proof
An element $x_0$ is a weak minimizer  of $f$, iff $f(x_0)\oplus U\nsubseteq f(x)$ is satisfied for all $U\in \mathcal U_Z$ and all $x\in X$. In other words, if and only if $0\notin \Int \of{f(x)\idif f(x_0)}$ is satisfied.
Obviously, if this is not satisfied, then there exists  $x\in X$ such that $0\in\Int \of{f(x)\idif f(x_0)}\subseteq \Int f'(x_0,x-x_0)$. Hence, if $x_0$ solves the variational inequality, then $x_0$ is a weak minimizer of $f$.
On the other hand, if $x_0$ is a weak minimizer of $f$, then especially for all $x\in X$ and all $t>0$ it holds
$0\notin\Int\frac{1}{t}\of{f(x_0+t(x-x_0))\idif f(x_0)}$, hence by Lemma \ref{prop:DirDer_of_conv_SV}
\[
0\notin\bigcup\limits_{t>0}\Int\frac{1}{t}\of{f(x_0+t(x-x_0))\idif f(x_0)}=\Int f'(x_0,x-x_0).
\]
\pend

In Section \ref{sec:Setting} we introduced also a scalarization of \eqref{eq:WMin}, thorough condition \eqref{eq:ScWMin}. The following results proves that, under some regularity condition, we have also equivalence between the scalarized optimization problem and variational inequalities.

\begin{theorem}\label{thm:Stamp_Min2}
Let $f:X\to G^\triup$ be a convex function, $x_0\in \dom f$.
If $x_0$ solves the scalarized Stampacchia variational inequality \eqref{eq:ScWSVI}, then it satisfies \eqref{eq:ScWMin}.
\end{theorem}
\proof
Since each scalarization $\vp_{f,z^*}:X\to\OLR$ is convex, $0\leq \vp'_{f,z^*}(x_0,x-x_0)$ implies $\vp_{f,z^*}(x_0)\leq \vp_{f,z^*}(x)\neq-\infty$. Hence if $x_0$ solves the Stampacchia variational inequality \eqref{eq:ScWSVI}, then for all $x\in X$ there exists a $z^*\in W^*$ such that $\vp_{f,z^*}(x_0)\leq \vp_{f,z^*}(x)\neq-\infty$ is satisfied and therefore $x_0$ satisfies \eqref{eq:ScWMin}.

\pend

The reverse implication needs further assumptions to hold.

\begin{theorem}\label{thm:Stamp_Min3}
Let $f:X\to G^\triup$ be a convex function, $x_0\in \dom f$.
If $x_0$ satisfies \eqref{eq:ScWMin} and any of the following conditions is satisfied:
\begin{enumerate}[(a)]
\item The regularity assumption \eqref{eq:reg_sharp} is satisfied;
\item It exists a finite subset $M^*\subseteq W^*$ such that
\[
\forall x\in X\, \exists z^*\in M^*:\quad \vp_{f,z^*}(x_0)\leq \vp_{f,z^*}(x)\neq-\infty;
\]

%\item  $f_{x_0,x}$ is lower $W^*$--equicontinuous in $0$ for all $x\in X$ and $\vp_{f,z^*}(x_0)\in\R$ is satisfied for all $z^*\in W^*\cap(0^+f(x_0))^-$.

\end{enumerate}
then $x_0$ solves \eqref{eq:ScWSVI}
\end{theorem}
\proof
\begin{enumerate}[(a)]

\item If $x_0$ satisfies \eqref{eq:ScWMin} and $f(x_0)\neq Z$, then it satisfies \eqref{eq:WMin} and, by Theorem \ref{thm:Stamp_Min1}, this implies that $x_0$ solves the Stampacchia variational inequality \eqref{eq:WSVI}.
If additionally the regularity assumption \eqref{eq:reg_sharp} is satisfied, then by Lemma \ref{lem:SVI} this implies that $x_0$ solves \eqref{eq:ScWSVI}.

\item
Let $x\in X$ be given. Then for all $t\in\of{0,1}$ there exists a $z^*\in M^*$ such that $\vp_{f,z^*}(x_0)\leq \vp_{f,z^*}(x_0+t(x-x_0))\neq-\infty$. As $M^*$ is finite, there exists a $z_0^*\in M^*$ and a sequence $t_n\downarrow 0$ in $\of{0,1}$ such that
$\vp_{f,z_0^*}(x_0)\leq \vp_{f,z_0^*}(x_0+t_n(x-x_0))\neq-\infty$, hence by convexity either $\sqb{x_0,x}\cap\dom f=\cb{x_0}$ and $\vp'_{f,z^*}(x_0,x-x_0)=+\infty$, or $\vp_{f,z^*_0}(x_0)\neq-\infty$, and
\[
0\leq \inf\limits_{n\in \N}\frac{1}{t_n}\of{\vp_{f,z^*}(x_0+t_n(x-x_0))\idif \vp_{f,z^*}(x_0)}
\]
and as $t_n$ converges to $0$, this implies $0\leq \vp'_{f,z^*}(x_0,x-x_0)$, hence
 $x_0$ solves \eqref{eq:ScWSVI}.

\end{enumerate}

\pend

%%%%%%%%%%%%%%%%%%%%%%%%%%%%%%%%%%

The study of variational inequalities related to optimization problems is classically divided into two parts. The first one relates to Stampacchia-type inequalities (see eg. \cite{Giannessi80}) and the second to Minty-type (see eg. \cite{Giannessi97}). Indeed, the differentiable Minty--type variational inequality, roughly speaking, evaluates the directional derivatives at some point $x$ along the direction $u=x_0-x$. This motivates the following definition.

\begin{definition}\label{def:MVI}
Let $f:X\to\G^\triup$ be a convex function and $f':X\times X\to\G^\triup$ its directional derivative.
Then $x_0$ is said to be a solution to the {\em weak Minty variational inequality}, iff $x_0\in\dom f$ and
\begin{align}\tag{W-MVI}\label{eq:WMVI}
f(x_0)=Z\quad\vee\quad \forall x\in X:\quad f'(x,x_0-x)\nsubseteq \Int 0^+f(x).
\end{align}
\end{definition}

As for Definition \ref{def:SVI}, we can provide a scalarization of \eqref{eq:WMVI} in the following lemma. However, a complete equivalence holds only for set-valued extensions of convex vector-valued function.

\begin{lemma}\label{lem:MVI}
If $x_0\in\dom f$ satisfies property \eqref{eq:WMVI}, then it also satisfies
\begin{align}\tag{Sc-W-MVI}\label{eq:ScWMVI}
f(x_0)=Z\quad\vee\quad \forall x\in X\, \exists z^*\in W^*:\quad \vp_{f,z^*}(x)\neq-\infty\wedge \vp'_{f,z^*}(x,x_0-x)\leq 0.
\end{align}
Moreover, let $f(x)=\psi^C(x)$ be true for some convex $\psi:S\subseteq X\to Z$ for all $x\in X$.
If $C$ is either Daniell, or $C\cap(k+(-C))$ is compact for all $k\in \Int C$, then
then for all $x\in S$ it exists $k_0\in \Int C$ such that $f'(x,x_0-x)=k_0+C$ and equivalence holds true.
\end{lemma}
\proof
By a separation argument, if $f'(x,x_0-x)\nsubseteq \Int 0^+f(x)$ is satisfied then either $f(x)=\emptyset$ and $f'(x,u)=Z$, in which case the statement is satisfied, or there exist $z^*\in W^*$ and $z \in f'(x,x_0-x)$ such that $\vp_{f,z^*}(x)\neq-\infty$ and $-\sigma(z^*|f'(x,x_0-x))\leq -z^*(z)\leq 0$, compare Lemma \ref{lem:Rec_A}.
By Proposition \ref{prop:vp'fLeqf'} the inequality $\vp'_{f,z^*}(x,x_0-x)\leq -\sigma(z^*|f'(x,x_0-x))$ is always satisfied, hence \eqref{eq:WMVI} implies \eqref{eq:ScWMVI}.

On the other hand if $f$ is the set--valued extension of a convex vector--valued function, then by Proposition \ref{prop:vp'fLeqf'2} $\vp'_{f,z^*}(x,x_0-x)= -\sigma(z^*|f'(x,x_0-x))$ is true for all $z^*\in W^*$ and applying Lemma \ref{lem:Rec_A} and Proposition \ref{prop:A_nsubset_intB} proves that \eqref{eq:ScWMVI} implies
\begin{align}
f(x_0)=Z\quad\vee\quad \forall x\in X\,\forall U\in\mathcal U_Z:\; f'(x,x_0-x)\oplus U\nsubseteq 0^+f(x).
\end{align}
It is left to prove that under the given assumptions this implies \eqref{eq:WMVI}.
If $x\notin \dom f$ or $0^+f'(x,x_0-x)\neq C$, then there is nothing to prove. Hence, let $x\in S$ and $0^+f'(x,x_0-x)=C$.\\
If $f'(x,x_0-x)\subseteq \Int C$, then $k_t=\frac{1}{t}\of{\psi(x+t(x_0-x))-\psi(x)}$ is a monotonly decreasing net in $\Int C$ as $t$ converges towards $0$ and bounded from below by $0\in Z$.

If $C$ is Daniell, this implies the differential quotient converges towards $k_0=\inf\limits_{t>0}k_t\in C$ and $f'(x,x_0-x)= k_0+ C$. Hence especially there exists a neighbourhood $U\in\mathcal U(0)$ with $k_0+U\subseteq C$, proving the equivalence.

If $C\cap(k+(-C))$ is compact for all $k\in\Int C$, then this is especially true for all $k_t$, $t>0$ and there exists a convergent subnet $k_{t_i}\to k_0$ with $k_0\in C$ and $k_0$ is a lower bound of $\cb{k_t}_{t>0}$. Hence, $f'(x,x_0-x)\subseteq k_0+C$ and $k_0\in\cl\bigcup\limits_{t>0}k_t+C$ proves $f'(x,x_0-x)=k_0+C\subseteq \Int C$. Hence especially there exists a neighbourhood $U\in\mathcal U(0)$ with $k_0+U\subseteq C$, proving the equivalence.

%If the set of differential quotients possesses a convergent subnet with limit $z_0$, then $z_0$ is a lower bound of the set $f'(x,x_0-x)$ and $z_0\in\Int C$, hence there exists a neighborhood $U\in\mathcal U_Z$ such that
%\[
%f'(x,x_0-x)\oplus U\subseteq \cl(\cb{z_0}+U)\subseteq C.
%\]
%Hence let $f'(x,x_0-x)\subseteq \Int C$ be assumed and the set $\cb{\frac{1}{t}\of{\psi(x+t(x_0-x))-\psi(x)}| t>0}$ does not possess any convergent subnets.
%
%For all $t>0$, define
%\[
%U_t=\of{C\cap \of{\frac{1}{t}\of{\psi(x+t(x_0-x))-\psi(x)}+(-C)}}+\cb{-\frac{1}{2t}\of{\psi(x+t(x_0-x))-\psi(x)}},
%\]
%a neighborhood of $0$. It holds $U_s\subseteq U_t$  whenever $0<s<t$ and $0\in \bigcap\limits_{0<t} U_t=U_0$.
%First, assume $\Int U_0=\emptyset$. Then for all $U\in\mathcal U_Z$ there is a $t>0$ such that
%$\cb{\frac{1}{2t}\of{\psi(x+t(x_0-x))-\psi(x)}}+U\nsubseteq C$, hence the set
%$\cb{\frac{1}{2t}\of{\psi(x+t(x_0-x))-\psi(x)}| t>0}$  possess a convergent subnet with limit outside of $\Int C$, a contradiction.
%Therefore, $U_0\in\mathcal U_Z$ is proven and
%\[
%\cb{\frac{1}{t}\of{\psi(x+t(x_0-x))-\psi(x)}}+U_0\subseteq C
%\]
%for all $0<t$. Hence,
%\[
%f'(x,x_0-x)\oplus U_0=\cl\bigcup\limits_{0<t}\of{\cb{\frac{1}{t}\of{\psi(x+t(x_0-x))-\psi(x)}}+C+U_0}\subseteq C.
%\]

\pend

Notably, as $C$ is closed, if $Z$ is finite dimensional, then $C\cap (k+(-C))$ is closed and bounded, hence compact for all $k\in Z$.

In the general setting of problem \eqref{P}, to prove the variational inequality characterization of weak minimizers, we need to apply a scalarization argument. Therefore we begin to study the scalarized version of the Minty inequality. Indeed, the next propositions show that the solution set to \eqref{eq:ScWMin} is always a subset of the solutions of \eqref{eq:ScWMVI}, while equality is satisfied under additional regularity assumptions.

\begin{theorem}\label{thm:Minty_Min1}
Let $f:X\to G^\triup$ be a convex function, $x_0\in \dom f$.
 If $x_0$ satisfies \eqref{eq:ScWMin} then it solves \eqref{eq:ScWMVI}.
\end{theorem}
\proof
If $x_0$ satisfies \eqref{eq:ScWMin} then either $f(x_0)=Z$, or for all $x\in X$ there exists a $z^*\in W^*$ such that $\vp_{f,z^*}(x)\neq-\infty$ and
\[
\vp'_{f,z^*}(x,x_0-x)\leq \vp_{f,z^*}(x_0)\idif \vp_{f,z^*}(x)\leq 0.
\]
\pend

\begin{theorem}\label{thm:Minty_Min2}
Let $f:X\to G^\triup$ be a convex function and $x_0\in \dom f$ solves \eqref{eq:ScWMVI}. If it exists a finite subset $M^*\subseteq W^*$ such that $f_{x_0,x}$ is $M^*$-l.s.c. in $0\in\dom f_{x_0,x}$ for all $x\in X$ and %either of the following two regualrity assumptions are satisfied:
%\begin{enumerate}[(a)]
%\item

\[
\forall x\in X\, \exists z^*\in M^*:\quad \vp_{f,z^*}(x)\neq-\infty\wedge \vp'_{f,z^*}(x,x_0-x)\leq 0;
\]
%
%\item  For all $z^*\in W^*\cap (0^+f(x_0))^-$ it holds $\vp_{f,z^*}(x_0)\in\R$ is true and $f_{x_0,x}$ is lower $W^*$-equicontinuous in $0\in\dom f_{x_0,x}$ for all $x\in X$
%
%
%\end{enumerate}
then $x_0$ satisfies \eqref{eq:ScWMin}.
\end{theorem}
\proof
%\begin{enumerate}[(a)]
%\item
Let $x\in X$ be given and $x_t=x_0+t(x-x_0)$. By convexity of $f$, if $\vp_{f,z^*}(x_t)\neq-\infty$ and $\vp'_{f,z^*}(x_t,x_0-x_t)\leq 0$, then  $\vp_{f,z^*}(x)\neq-\infty$ and $\vp'_{f,z^*}(x,x_0-x)\leq 0$ is satisfied and $\vp_{f,z^*}(x_t)\leq \vp_{f,z^*}(x)$.
As by assumption the set $M^*$ is finite, for any $x\in X$ there exists a $z^*\in M^*$ such that for all $t>0$ it holds $\vp_{f,z^*}(x_t)\neq-\infty$ and $\vp'_{f,z^*}(x_t,x_0-x_t)\leq 0$.
As $(\vp_{f,z^*})_{x_0,x}$ is convex and l.s.c. in $0$, this implies $\vp_{f,z^*}(x_0)=\inf\limits_{t\in\sqb{0,1}}\vp_{f,z^*}(x_t)\leq \vp_{f,z^*}(x)$.

%\item
%If either $f(x_0)=Z$ or $x\notin\dom f$, then obviously \eqref{eq:ScWMin} is satisfied in $x$. Hence let $x\in\dom f$ and $f(x_0)\neq Z$. For any $0\leq t$, we set $x_t=x_0+t(x-x_0)$. By assumption, for any $0<t$ there exists a $z^*_t\in W^*$ such that
%$\vp_{f,z^*_t}(x_t)\neq-\infty$ and $\vp'_{f,z^*_t}(x_t,x_0-x_t)\leq 0$. Especially, as $f_{x_0,x}$ is convex and l.s.c. on the interval $\sqb{0,t}$, the infimum $-\sigma\of{z^*_t| f\sqb{x_0,x}}\neq-\infty$ of $\vp_{f,z^*_t}$ is attained at some $\hat x=(h_tx_t+(1-h_t)x_0)$ with $h_t\in\sqb{0,1}$ and for all $t\leq s$ it holds $\vp_{f,z^*_t}(x_s)\neq-\infty$, $\vp'_{f,z^*_t}(x_s,x_0-x_s)\leq  0$.
%
%If $\vp_{f,z_t^*}(x_0)=\vp_{f,z^*_t}(h_tx_t+(1-h_t)x_0)$ then we are finished, hence assume $\vp_{f,z^*_t}(h_tx_t+(1-h_t)x_0)<\vp_{f,z_t^*}(x_0)$.
%
%By the same arguments as in the proof of proposition \ref{thm:Stamp_Min3} (c), we can show that
% for any $x\in X$ there exists a $z_0^*\in W^*\cap (0^+f(x_0))^-$ such that $\vp_{f,z_0^*}$ attains its finite infimum over the intervall $\sqb{x_0,x}$ at $x_0$, hence $\vp_{f,z_0^*}(x_0)\leq \vp_{f,z_0^*}(x)\neq-\infty$.
%Therefore $x_0$ satisfies \eqref{eq:ScWMin}.
\pend
%\end{enumerate}

Based on the previous results we can prove equivalence between solutions of Minty type inequality and weak minimizers at least when $f=\psi^C$.

\begin{corollary}\label{cor:MVI_Min}
Let $f:X\to G^\triup$ be a convex function.
\begin{enumerate}[(a)]
\item
If $f$ satisfies the regualarity assumption given in Theorem \ref{thm:Minty_Min2} for $x_0\in\dom f$ and $x_0$ solves \eqref{eq:WMVI}, then it also satisfies \eqref{eq:WMin}.
\item
If  $Z$ has finite dimension or $C$ is Daniell, then if $f(x)=\psi^C(x)$ for all $x\in X$ and $x_0\in\dom f$ satisfies \eqref{eq:WMin}, then $x_0$ also solves \eqref{eq:WMVI}.
\end{enumerate}
\end{corollary}
\proof
The implication in $(a)$ is proven in Lemma \ref{lem:MVI}, Theorem \ref{thm:Minty_Min2} and Proposition \ref{prop:A_nsubset_intB}, while the implication in $(b)$ is a corollary of Proposition \ref{rem:wEff_wMin}, Theorem \ref{thm:Minty_Min1} and Lemma \ref{lem:MVI}.
\pend

Another regularity assumption than that in Theorem \ref{thm:Minty_Min2} can be found in \cite{CrespiRoccaSchrage14}.

The following example shows that we cannot obtain a result similar to Theorem \ref{thm:Stamp_Min1} for Minty type variational inequality.

\begin{example}
Consider $Z=\R^2$, ordered by the natural ordering cone $C=\R^2_+$ and $X=\R$. The function $f:X\to\G^\triup$  given by
\[
f(t)=\begin{cases}
\cb{(z_1,z_2)^T\in Z\st -t\leq z_1,z_2,\, t\leq z_1+z_2},&\text{ if }t\in\of{0,1};\\
\emptyset,&\text{ elsewhere}
\end{cases}
\]
is  convex and $\C$--l.s.c. everywhere.
Then $f'(1,-1)=(1,1)^T+C$, hence there exists  $t\in\dom f$ and $U\in\mathcal U_Z$ such that $f'(t,0-t)+U\subseteq 0^+f(t)=C$ and obviously $f'(t,0-t)\subseteq\Int 0^+f(t)$. However, $f(0)\nsubseteq \Int f(t)$ for all $t\in\R$, hence $f(0)$ is a weak-l-minimal element of $f\sqb{X}$ and thus especially satisfies \eqref{eq:WMin} and \eqref{eq:ScWMin}, but the Minty variational inequality \eqref{eq:WMVI} is not satisfied.
\end{example}

To summarize, we have proved the following chain of characterization of weak minimizers of problem \eqref{P} for convex functions through set-valued variational inequalities.

	\begin{center}
		\includegraphics[width=14cm]{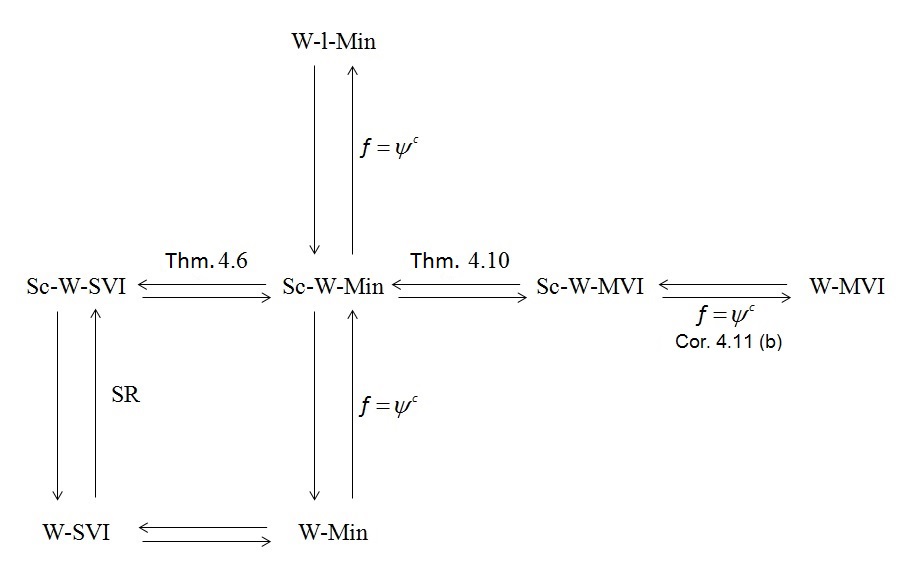}
	\end{center}

Finally, when additionally $f=\psi^C$ we can simplify the previous results to gain a characterization of weak efficiency in vector optimization.\\
In this case, the scalarized variational inequality \eqref{eq:ScWSVI} is equivalent to its set--valued counterpart without further assumptions. If the ordering cone is Daniell, then the same is true for \eqref{eq:ScWMVI}  and \eqref{eq:WMVI}.

\begin{proposition}
Let $\psi:S\subseteq X\to Z$ be a $C$--convex function, $x_0\in S$ and $f(x)=\psi^C(x)$ for all $x\in X$. Then
\begin{enumerate}[(a)]
\item the Stampacchia variational inequalities of type \eqref{eq:WSVI} and \eqref{eq:ScWSVI} are equivalent;
\item the Minty variational inequalities of type \eqref{eq:WMVI} and \eqref{eq:ScWMVI} are equivalent, if $Z$ has finite dimension or $C$ is Daniell.
\end{enumerate}
\end{proposition}
\proof
\begin{enumerate}[(a)]
\item Assuming $f(x)=\psi^C(x)$ for all $x\in X$ is true, the regularity assumption \eqref{eq:reg_sharp} is satisfied and equivalence follows from  Lemma \ref{lem:SVI}.
\item This is Lemma \ref{lem:MVI}.
\pend
\end{enumerate}

Finally, we provide the classical chain of relations for weak efficiency (compare e.g. \cite{Giannessi97, CGR}) for $C$--convex functions as corollaries of the results proved in the general case.

	\begin{center}
		\includegraphics[width=10cm]{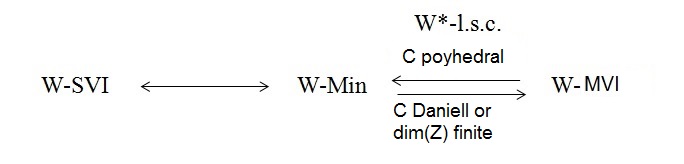}
	\end{center}

The following corollaries state the implications in the scheme.

\begin{corollary}\label{cor:Stamp_wEff}
Let $\psi:S\subseteq X\to Z$ be a $C$--convex function and $f(x)=\psi^C(x)$ for all $x\in X$. Then
$x_0\in S$ solves \eqref{eq:WSVI} if and only if $x_0$ is a weakly efficient solution of the vector optimization problem \eqref{VOP}.
\end{corollary}

\begin{corollary}\label{cor:Minty_wEff}
Let $\psi:S\subseteq X\to Z$ be a $C$--convex function and $f(x)=\psi^C(x)$ for all $x\in X$.
\begin{enumerate}[(a)]
\item
If $x_0\in S$ is a weakly efficient solution of the vector optimization problem \eqref{VOP},
then $x_0$ solves \eqref{eq:ScWMVI};
\item
If additionally $f_{x_0,x}$ is $\C$--l.s.c in $0$ for all $x\in X$, $C$ is polyhedral, then $x_0\in S$ solves \eqref{eq:ScWMVI} if and only if $x_0$ is a weakly efficient solution of the vector optimization problem \eqref{VOP}.
\end{enumerate}
\end{corollary}
\proof
\begin{enumerate}[(a)]
\item
If $x_0\in S$ is a weakly efficient solution of the vector optimization problem \eqref{VOP},
then $x_0$ is a weak minimizer of $f$ this implies \eqref{eq:WMVI}
\item
The reverse implication follows from Theorem \ref{thm:Minty_Min2}, as $\vp_{f,z^*}(x_0)=-z^*(\psi(x_0))\in\R$ is true.
\pend
\end{enumerate}

If additionally either $Z$ has finite dimension or $C$ is Daniell, then the weak scalarized Minty variational inequality can be replaced by \eqref{eq:WMVI}, compare Lemma \ref{lem:MVI}.

The main advantage of these results, compared with those in \cite{CrespiGinRoc2005, Ginchev2007} is that  $\psi(x_0)\in\wEff \psi\sqb{X}$ is characterized using a Minty or Stampacchia type variational inequality for the epigraphical extension of $\psi$ saving us the effort of introducing ''infinite elements'' of $Z$ to cope with possible unboundedness of the differential quotient $\frac{1}{t}\of{\psi(x_0+tu)-\psi(x_0)}$.

%\section{Acknowledgment}

%We thank the anonymous referees for their valuable comments that helped to improve the manuscript.

\addcontentsline{toc}{section}{Bibliography}
%\bibliography{../../../../Literatur}

\begin{thebibliography}{10}

%\bibitem{AliprantisBorder}
%C.D. Aliprantis and K.C. Border.
%\newblock {\em {Infinite Dimensional Analysis: A Hitchhiker's Guide}}, Volume~4
%  of {\em Studies in Economic Theory}.
%\newblock Springer Verlag, 2006.

\bibitem{AnsariLee2010}
Q.H. Ansari and G.M. Lee.
\newblock Nonsmooth vector optimization problems and minty vector variational
  inequalities.
\newblock {\em Journal of Optimization Theory and Applications}, 145:1--16, 2010.

\bibitem{Aubin71}
	J.-P. Aubin.
 \newblock  	{\em A Pareto Minimal Principle} in
\newblock H.W. Kuhn and G.P. Szeg{\"o} (eds.).
\newblock 	Differential games and related topics,
\newblock North-Holland Pub. Co., 147--175, 1971.

\bibitem{AubinFrankowska}
J.-P. Aubin and H.~Frankowska.
\newblock {\em {Set-Valued Analysis}}, Volume~2 of {\em Systems \& Control:
  Foundations \& Applications}.
\newblock Birkh\"auser Boston Inc., Boston, MA, 1990.

\bibitem{cgr-erice}
G.P. Crespi, and I. Ginchev, I. and M. Rocca.
 {\sl Variational inequalities in vector optimization.} In: F. Giannessi and A. Maugeri eds.,  Variational Analysis and Applications,  Kluwer Acad. Publ., Dordrecht, 2004.

\bibitem{cgrJOTA}
G.P. Crespi, I.~Ginchev, and M.~Rocca.
\newblock Minty Variational Inequalities, Increase-Along-Rays Property and Optimization.
\newblock {\em Journal of optimization theory and applications}, 123(3): 479-496, 2004.

\bibitem{CrespiGinRoc2005}
G.P. Crespi, I.~Ginchev, and M.~Rocca.
\newblock First--order optimality conditions in set--valued optimization.
\newblock {\em Mathematical Methods of Operations Research}, 63(1):87--106,  2006.

\bibitem{CGR}
G.P. Crespi, I.~Ginchev, and M.~Rocca.
\newblock Some remarks on the Minty vector variational principle
\newblock {\em Journal of Mathematical Analysis and Applications}, 345(1): 165-175, 2008


\bibitem{CrespiHamelSchrage2015}
G.P. Crespi, A.H..~Hamel, and C.~Schrage.
\newblock A Minty variational principle for set optimization
\newblock {\em Journal of Mathematical Analysis and Applications}, 423(1): 770--796, 2015.


\bibitem{CrespiRoccaSchrage14}
 G.P. Crespi, M. Rocca and C. Schrage.
\newblock    Minty variational inequality and weak minimality in set optimization
\newblock {\em Journal of Optimization Theory and Applications}, 166(3): 804--24,	2015.
%DOI: 10.1007/s10957-014-0679-3

\bibitem{CrespiSchrage13a}
G.P. Crespi and C.~Schrage.
\newblock Set optimization meets variational inequalities.
in {\sl Set Optimization and Applications - The State of the Art} Hamel, A.H., Heyde, F., L\"ohne, A., Rudloff, B. and Schrage, C. (eds.), pages 213--247 Springer New York, 2015.


\bibitem{Fuchs66}
Fuchs, L.
\newblock {\em {Teilweise geordnete algebraische Strukturen}}.
\newblock Vandenhoeck u. Ruprecht, G\"ottingen, 1966.

\bibitem{FuchssteinerLusky1981}
    Fuchssteiner, B. and Lusky, W.
    \newblock{\em{Convex Cones, volume 56 of Mathematics Studies}}
    \newblock Elsevier Science Ltd, Amsterdam, 1981



\bibitem{galatos2007residuated}
N.~Galatos.
\newblock {\em {Residuated Lattices: An Algebraic Glimpse at Substructural
  Logics}}.
\newblock Elsevier Science Ltd, Amsterdam, 2007.

\bibitem{GetanMaLeSi}
J.~Getan, J.~E. Martinez-Legaz, and I.~Singer.
\newblock {(*, s)-Dualities}.
\newblock {\em Journal of Mathematical Sciences}, 115(4):2506--2541, 2003.

\bibitem{Ginchev2007}
I.~Ginchev.
\newblock Vector optimization problems with quasiconvex constraints.
\newblock {\em Journal of Gobal Optimization}, 44:111--130, 2007.

\bibitem{Giannessi80}
F. Giannessi.
 Theorems of the alternative, quadratic programs and complementarity problems, in {\sl Variational Inequalities and Complementarity Problems. Theory and applications} R.W.~Cottle, F.~Giannessi, J.L.~Lions (eds.), Wiley, New York, pp. 151-186, 1980.

\bibitem{Giannessi97}
F. Giannessi.
  {\sl On Minty variational principle.} New trends in mathematical programming, Kluwer Acad. Publ., Dordrecht, 93-99, 1997.

\bibitem{giorgi1992dini}
G~Giorgi and S~Koml{\'o}si.
\newblock Dini derivatives in optimization --Part i.
\newblock {\em Decisions in Economics and Finance}, 15(1):3--30, 1992.

%\bibitem{GRTZ}
%A.~G{\"o}pfert, H.~Riahi, Ch. Tammer, and C.~Z{\u{a}}linescu.
%\newblock {\em {Variational Methods in Partially Ordered Spaces}}.
%\newblock CMS Books in Mathematics/Ouvrages de Math\'ematiques de la SMC, 17.
 % Springer-Verlag, New York, 2003.



\bibitem{Hamel09}
A.~H. Hamel.
\newblock {A Duality Theory for Set-Valued Functions I: Fenchel Conjugation
  Theory}.
\newblock {\em Set-Valued and Variational Analysis}, 17(2):153--182, 2009.

\bibitem{HamelHabil05}
A.~H. Hamel.
\newblock Variational {P}rinciples on {M}etric and {U}niform {S}paces.
\newblock Habilitationsschrift, {Halle}, 2005.

\bibitem{FivePersonSoMF}
Hamel, A.H., Heyde, F., L\"ohne, A., Rudloff, B. and Schrage, C.
\newblock Set OptimizationA Rather Short Introduction in {\sl Set Optimization and Applications - The State of the Art} Hamel, A.H., Heyde, F., L\"ohne, A., Rudloff, B. and Schrage, C. (eds.), pages 65--141 Springer New York, 2015.



\bibitem{HamelSchrage12}
A.~H. Hamel and C.~Schrage.
\newblock Notes on extended real- and set-valued functions.
\newblock {\em Journal of Convex Analysis}, 19(2):355--384, 2012.

\bibitem{HamelSchrage13PJO}
A.~H. Hamel and C.~Schrage.
\newblock Directional derivatives, subdifferentials and optimality
conditions for set-valued convex functions.
\newblock {\em Pacific Journal of Optimization}, 10(4): 667-684, 2014.

\bibitem{HernandezMarin05}
E.~Hern{\'a}ndez and L.~Rodr{\'{\i}}guez-Mar{\'{\i}}n.
\newblock Nonconvex scalarization in set optimization with set-valued maps.
\newblock {\em Journal of Mathematical Analysis and Applications}, 325(1):1--18, 2007.

\bibitem{HeydeLoehne11}
F.~Heyde and A.~L\"ohne.
\newblock Solution concepts in vector optimization. a fresh look at an old
  story.
\newblock {\em Optimization}, 60(12):1421--1440, 2011.

\bibitem{HeydeSchrage11R}
F.~Heyde and C.~Schrage.
\newblock Continuity of set--valued maps and a fundamental duality formula for
  set--valued optimization.
\newblock {\em Journal of Mathematical Analysis and Applications},
  397(2):772--784, 2013.

\bibitem{Kuroiwa98-1}
D.~Kuroiwa.
\newblock {The natural criteria in set--valued optimization}.
%\newblock{Investigations on nonlinear analysis and convex analysis (Japanese), Kyoto}.
\newblock{ RIMS Kokyuroku} 1031,  85--90, 1998.

\bibitem{kutateladze1972minkowski}
Kutateladze, S.S. and Rubinov, A.M..
\newblock{Minkowski duality and its applications}.
\newblock{Russian Mathematical Surveys} 27 (3), 137191, 1972

\bibitem{Loehne11Book}
A.~L\"ohne.
\newblock {\em {Vector Optimization with Infimum and Supremum}}.
\newblock Springer-Verlag, Berlin, 2011.

\bibitem{Luc89}
D. T. Luc.
\newblock {\em Theory of Vector Optimization}.
\newblock Lecture notes in Economics and Mathematical Systems, Springer-Verlag Berlin, 1989.


\bibitem{MartinezLegazSinger95}
J.E. Mart{\i}nez-Legaz and I.~Singer.
\newblock {Dualities associated to binary operations on $\bar R$}.
\newblock {\em Journal of Convex Analysis}, 2(1-2):185--209, 1995.

%\bibitem{Rockafellar}
%   R.T. Rockafellar.
%	\newblock {\em Convex Analysis}.
% 	\newblock Princeton Mathematical Series, No. 28, Princeton University Press, Princeton, N.J., 1970.

%\bibitem{RockafellarWets98}
%R.T. Rockafellar and R.J.B. Wets.
%\newblock {\em {Variational Analysis}}.
%\newblock Springer, 2004.

\bibitem{Diss}
Schrage, C.
\newblock{Set-Valued Convex Analyis},
\newblock  PhD-THesis, Halle, 2009.

\bibitem{Schrage10Opt}
Schrage, C.
\newblock Scalar representation and conjugation of set--valued functions.
\newblock {\em Optimization}, 64(2): 197--223, 2015.


\bibitem{YYT}
X.M. Yang and X.Q.Yang and K.L. Teo.
  {\sl Some remarks on the Minty vector variational inequality.} Journal of Optimization Theory and Applications, 121, 193-201, 2004.

\bibitem{Zalinescu02}
C.~Z{\u{a}}linescu.
\newblock {\em {Convex Analysis in General Vector Spaces}}.
\newblock World Scientific Publishing Co. Inc., River Edge, NJ, 2002.

\end{thebibliography}

\end{document}